\documentclass{siamltex}

\usepackage{amsmath}
\usepackage{amssymb}

\usepackage[toc,page]{appendix}

\usepackage{tikz}
\usepackage{pgfplots}
\usepgfplotslibrary{groupplots}
\usepackage[ComplainAboutUnexternalized]{externaltikz}
\usepackage{calc}

\usepackage{subcaption}
\usepackage{hyperref}
 \usepackage{mathtools}

\newtheorem{remark}{Remark}

\newcommand{\rhoB}{\rho_B} 
\newcommand{\rhoK}{\rho_K}

\begin{document}
	
\title{A kinetic traffic network  model and its macroscopic limit: merging lanes}

\author{R. Borsche\footnotemark[1] 
	\and  A. Klar\footnotemark[1] \footnotemark[2]}
\footnotetext[1]{Technische Universit\"at Kaiserslautern, Department of Mathematics, Erwin-Schr\"odinger-Stra{\ss}e, 67663 Kaiserslautern, Germany 
	(\{borsche, klar\}@mathematik.uni-kl.de)}
\footnotetext[2]{Fraunhofer ITWM, Fraunhoferplatz 1, 67663 Kaiserslautern, Germany}

	\maketitle

	\begin{abstract}
		In this paper we propose coupling conditions for a kinetic two velocity model for vehicular traffic on networks. These conditions are based on the 
		consideration of the free space on the respective roads. The macroscopic limit of the kinetic relaxation system is a  classical scalar conservation law for traffic flow.
		Similar to the asymptotic limit of boundary value problems for kinetic models, we consider here the limit of the full network problem including the coupling conditions at the nodes.
	 	An asymptotic analysis of the interface  layers at the nodes and  a matching procedure	using half-Riemann problems for the limit conservation law are used to derive coupling conditions for classical macroscopic traffic models on the network from the kinetic ones.	
	\end{abstract}

\section{Introduction}

Modeling and simulation of traffic flow on road networks, has been investigated intensively using hyperbolic partial differential equations. 
Different models have been used, ranging from scalar conservation laws like the Lighthill Whitham Richards model, through 
models using system of conservation laws, to kinetic descriptions of the flow. 
Derivations of these models from the underlying models in such hierarchies have been discussed as well. See,  \cite{AR,AKMR,Ber,Hel,HPRV20,Dag1,Dag2,KW97,L,PSTV17,PSTV171} for a non-exhaustive list of references.
To obtain a model for the dynamics on the full network, all these models have to be supplemented with coupling conditions at the nodes of the network. 
Coupling conditions for scalar conservation laws and systems of conservation laws on networks  have been discussed in many papers, see, for example,  \cite{Cor17,BNR14,HR95,GPBook,LS02,CGP05,BHK06a,CG08,CM08,G10}.
Kinetic equations on networks have been considered, for example  in  \cite{HM09,BKKP16}. 
The present work aims at giving a derivation of a network model, based on a scalar conservation law with  coupling conditions at the nodes,  from an underlying kinetic or relaxation model with suitable coupling conditions at the nodes.

To explain the general procedure in  more detail, we consider a scaled kinetic problem or a general relaxation equation in 1D involving a scaling parameter $\epsilon$, which  converges for $\epsilon \rightarrow 0$ to  an associated scalar conservation law for traffic flow.
If such equations are considered on a network, it is sufficient to study a single coupling point or node, where  coupling conditions are required. 
Suitable  coupling conditions have  to be imposed for the kinetic problem at each node, see e.g. \cite{BKKP16}.
If on the network $\epsilon$ is send to zero, boundary layers at the junctions can arise.
To consider the limit  $\epsilon \rightarrow 0$ of the kinetic problem on the network, one has to proceed similarly
as in the case of kinetic boundary value problems, where  a complete picture of the convergence is only obtained, once so called boundary- and initial layers are investigated.
We refer to  \cite{BSS84,BLP79,G08,UTY03} for such a procedure in the case of kinetic equations and to  
\cite{WY99,WX99,LX96,X04} for the case of hyperbolic relaxation systems.  
However, the consideration of the kinetic layers alone is not sufficient to solve the coupling problem as waves can emerge from the junction:
the kinetic  half space problems have to be coupled to half Riemann problems for the macroscopic equation, as outlined in \cite{BK18b}.

Such a  procedure to derive coupling conditions for macroscopic equations from the underlying kinetic ones has been  discussed for linear systems in \cite{BK18c} using an asymptotic analysis of the situation near the nodes.
Using a detailed analysis of the kinetic boundary layers at a node, coupling conditions for the macroscopic equations have been  derived from kinetic coupling conditions in a simplified, but nonlinear case in \cite{BK18b}.

In the present work, we consider the case of a basic nonlinear kinetic traffic flow equation on a network with a small parameter $\epsilon$.
We propose  coupling conditions for this model and derive classical coupling conditions for nonlinear scalar traffic equations in the limit 
$\epsilon $ going to $0$. 
The investigation in this paper is restricted to the case of nodes with three edges with two ingoing and one outgoing road, i.e. a node where two lanes merge into one lane.
The case of diverging lanes will be discussed in a subsequent publication \cite{BK20}.

The paper is organized in the following way.  
In section \ref{equations} we present the kinetic traffic relaxation model, compare \cite{BK18}, and the associated scalar conservation law. 
In section \ref{sec:kineticcouplingconditions} coupling conditions for the kinetic relaxation model based on the available free space for cars entering a junction are discussed.
In section \ref{boundaryconditions} kinetic boundary layers are investigated, as well as the combination of the layer solutions with suitable Riemann solvers. 
This leads to boundary conditions for the scalar traffic equations depending on the kinetic boundary condition.
Then, the derivation of coupling conditions for the scalar hyperbolic limit problem from the kinetic coupling conditions is discussed.
Explicit coupling conditions for the macroscopic equation based on the kinetic coupling conditions are derived in section \ref{macroscopiccc}.	
The  derivation shows that a classical merge condition  for a nonlinear scalar conservation law can be interpreted on the kinetic level as a combination of the balance of fluxes and the equality of densites on all roads.
Finally, the solutions of the macroscopic equations on the network are compared numerically to the full solutions of the kinetic equation on the network in section \ref{Numerical results}.

\section{Kinetic and macroscopic traffic equations}
\label{equations}

We consider a minimal kinetic discrete velocity model \cite{BK18} with  two velocities $v_0 = 0$ and $v_1=1$.
$f_0$ represents the density of stopped cars, while $f_1$ is the density of driving cars. 
We define the total density of cars as $\rho = f_0+f_1\in [0,1]$ and the mean flux as $q =v_0 f_0 +v_1 f_1= f_1 $.
Reversely we obtain
$$
f_0 =\rho-q\ ,\quad  f_1 =q\ .
$$
The discrete velocity model is given by the following equations, see \cite{BK18}
\begin{eqnarray}
\label{eq:dvm01}
\begin{aligned}
\partial_t f_0  - \frac{1}{1-\rho }    f_1\partial_x f_0&= -\frac{1}{\epsilon} \left(f_0-\rho +F(\rho) \right)  \\
\partial_t f_1  + \partial_x f_1 + \frac{1}{1-\rho }   f_1\partial_x f_0&=  -\frac{1}{\epsilon} \left(f_1- F(\rho) \right)\ ,
\end{aligned}
\end{eqnarray}
where $F= F(\rho)$ is a given traffic density-flow function or fundamental diagram, i.e. a smooth function $F:[0,1]\rightarrow[0,1]$ with $F(0)=0=F(1)$ and $F^\prime(\rho) \le 1$ with graph  in
the triangle $0 \le \rho \le 1, 0 \le q \le \rho$.
In the following we restrict ourselves to strictly concave fundamental diagrams $F$ and denote
the point, where the maximum of $F$ is attained at the point $\rho^{\star}$ and the maximal value is 
$ F(\rho^{\star})=\sigma$.

The two eigenvalues corresponding to \eqref{eq:dvm01} are $\lambda_1 = - \frac{q}{1-\rho}\leq 0 <\lambda_2 = 1$.  
The respective eigenvectors  are
$r_1 = \left(
		1,
		\lambda_1 \right)^T, 
	r_2 = \left(
	1,1
	\right)\ $. 
A straightforward computation shows that the $r_1$- and  the $r_2$-field are both linearly degenerate. 
The integral curves (and shock curves) of the hyperbolic system are given by $q= q_L \frac{1-\rho}{1-\rho_L}$ for the 1-field and by $q = \rho -\rho_R+q_R$ for the 2-field. 
The region $0 \le \rho \le 1, 0 \le q \le \rho$ is an invariant region for the kinetic equations.
We refer to \cite{BK18} for details.

We note that the equations  can be rewritten in conservative form choosing the   variable 
$
	Z =  \frac{q}{1-\rho+q}\ .
$
Rewriting \eqref{eq:dvm01} we obtain
	\begin{align}\label{eq:lindeg+relax}
	\begin{aligned}
	\partial_t \rho + \partial_x q&=0\\
	\partial_t Z + \partial_x Z &= -\frac{(1-Z)}{\epsilon (1-\rho)} \left(q - F(\rho)  \right)\ 
	\end{aligned}
	\end{align}
with $q =  \frac{Z}{1-Z}(1-\rho)$.
A  Riemann invariant of the  second characteristic family is 
$$w =\rho - q = f_0 \in [0,1]\ .$$ 
A Riemann invariant of the first characteristic family is 
$$Z= \frac{q}{1-\rho+q} = \frac{q}{1-w}= \frac{f_1}{1-f_0}\in [0,1]\ .$$
We have $ 0 \le Z \le \rho$ and the fundamental diagram  $Z(\rho) = \frac{F(\rho)}{1-\rho+F(\rho)}$.

We interpret $1-w$ as the free space available or the maximal possible number of driving cars.
$Z$ can be understood as the ratio between the actual number of driving cars and the maximal possible number of driving cars.
	
Equation \eqref{eq:dvm01} can be also expressed in macroscopic variables $\rho$, $q$ as
\begin{align}\label{macro0}
\begin{aligned}
\partial_t \rho + \partial_x q &=0\\
\partial_t q + \frac{ q}{1-\rho}  \partial_x \rho  + (1-\frac{ q}{1-\rho})\partial_x  q  &=-\frac{1}{\epsilon} \left(q-F(\rho) \right) \ .
\end{aligned}
\end{align}
Concerning the convergence of its solutions towards the solutions of the scalar conservation law $\partial_t \rho + \partial_x F(\rho) =0$ as
$\epsilon $ tends to $0$ the subcharacteristic condition has to be satisfied \cite{LX96}.
Setting $q = F(\rho)$ in  the formula for the eigenvalues, the subcharacteristic condition states
$$
-  \frac{ F(\rho)}{1-\rho} \le F^\prime(\rho) \le 1 \  \mbox{ for } \  0 \le \rho \le 1\ .
$$

\begin{remark}
	The condition is fulfilled for strictly concave fundamental diagrams $F$.
	For example, in the classical LWR case with $F(\rho) =\rho (1-\rho) $ and  $F^\prime (\rho) = 1- 2 \rho$ the above  condition is
	$$
	-   \rho \le 1- 2 \rho \le 1\  \mbox{ for } \ 0 \le \rho \le 1\ ,
	$$
	which is obviously satisfied. 
\end{remark}

Finally we equip the kinetic problem \eqref{eq:dvm01} with boundary conditions.
At the left boundary at $x=x_L$ we have to prescribe a value for the 2- Riemann invariant $Z(x_L) = \frac{q(x_L)}{1-\rho(x_L)+q(x_L)} =  \frac{f_1(x_L)}{1-f_0(x_L)}$ and for the right boundary $x= x_R$ the 1-Riemann invariant $w(x_R)=f_0(x_R)$.
Note that the first eigenvalue is always non-positive and the second is constant, thus, the number of boundary conditions is fixed. 

A procedure to derive boundary conditions for the macroscopic equations from the kinetic ones is explained in section \ref{boundaryconditions}.

\section{Kinetic Coupling conditions} \label{sec:kineticcouplingconditions}
	
\begin{figure}[h!]
	\begin{center}
		\externaltikz{sketch_21node}{
			\begin{tikzpicture}[thick]
			\def\len{2}
			\node[fill,circle] (N) at (0,0){};
			\draw[->] (-\len,0.6)--(N) node[above,pos = 0.5]{$1$};
			\draw[->] (-\len,-0.6)--(N) node[below,pos = 0.5]{$2$};
			\draw[->] (0,0)--(\len,0) node[above,pos = 0.5]{$3$};
			\end{tikzpicture}
		}
	\end{center}
	\caption{A junction with two ingoing and one outgoing road (2-1 node). }
	\label{fig:merge}
\end{figure}
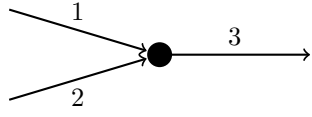

In this section we propose coupling conditions for the kinetic two-velocity model \eqref{eq:dvm01}.
They will be used to derive coupling conditions for the macroscopic equations in section \ref{macroscopiccc}.
We restrict ourselves here to the case of a junction with two ingoing and one outgoing lane, i.e. a junction with two merging lanes as shown in Figure \ref{fig:merge}.
The case of diverging lanes will be treated in a subsequent paper.
As on each road there is exactly one outgoing characteristic family, we have to provide  three conditions  at a junction connecting three roads. 
In any case the conservation of mass will be imposed, i.e. all cars entering a junction via one of the incoming roads
will exit on the outgoing road.
For the remaining two conditions we will analyze the maximal possible number of driving cars $1-w = 1-f_0$.  
Note that only the stopped cars $w = f_0$  on the outgoing roads can block the traffic, as the driving ones $f_1=q$ will clear the space.
Thus $1-f_0$ is the available free space for driving cars.
	
From the mathematical point of view it is important to supply values for the correct characteristic variables.
We denote with $\hat{\cdot}$ the known traces at the junction. The unknown characteristic variables as well as partially known quantities in the junctions do not have any indicator.
The roads are numbered as shown in Figure \ref{fig:merge}, the index of the road will be a superscript for the corresponding variables.

\subsection{Fair merging}
\label{sec:fair}	
We consider the free space on road $3$ for the cars entering the junction from road 1 and road 2.
For a car entering from road $1$, the free space on road $3$ is given by  $1-f_0^3=1-w^3$ minus the  cars entering from road $2$, i.e.  $f_1^2=q^2$.
Altogether we obtain	
\begin{align}
\label{eq:merge_kin_eq1}
1-f_0^1 = 1-f_0^3 - f_1^2\ .
\end{align}
Similarly, the cars entering from  road $2$ see the free space 
\begin{align*}
1-f_0^2 = 1-f_0^3 - f_1^1\ .
\end{align*}
As the mass should be conserved the third condition is 
\begin{align}\label{eq:merge_kin_eq3}
f_1^3 = q^3=  q^1+q^2= f_1^1+f_1^2\ .
\end{align}
These three equations (\ref{eq:merge_kin_eq1}-\ref{eq:merge_kin_eq3}) form a complete set of coupling conditions.
We rewrite these equations in the characteristic variables $w$ and $Z$
\begin{align*}
1-{w}^1&=(1-\hat{w}^3)-\hat{Z}^2 (1-{w}^2)\\
1-{w}^2&=(1-\hat{w}^3)-\hat{Z}^1(1-{w}^1)\\
{Z}^3(1-\hat{w}^3) &= \hat{Z}^1(1-{w}^1)+\hat{Z}^2(1-{w}^2)\ ,
\end{align*}
which allows to distinguish known from unknown quantities.
This is a linear system for the unknown variables $w^1,w^2,Z^3$. 
It is uniquely solvable for $\hat Z^1 \hat Z^2 \neq 1$ and $\hat{w}^3<1$ with the solution

\begin{align}
\label{z3}
1-w^{1/2}&= \frac{1-\hat Z^{2/1}}{1- \hat Z^{1} \hat Z^{2}} (1-\hat w^3) = \alpha^{1/2}(1-\hat w^3)  \nonumber \\
Z^3 &= \frac{\hat Z^{1}- 2 \hat Z^{1} \hat Z^{2}+ \hat Z^{2}}{1- \hat Z^{1} \hat Z^{2}}= \alpha_1 \hat Z^{1}+ \alpha_2 \hat Z^{2}\ ,
\end{align}	
where
$\alpha^{1/2}= \frac{1-\hat Z^{2/1}}{1- \hat Z^{1} \hat Z^{2}} \in [0,1]$.
In case $\hat{w}^3 = 1$ we have  $w^1=1=w^2$ and $Z^3$ is given again by \eqref{z3}.
In the particular case that $\hat Z^1=\hat Z^2=1$ the values of $\alpha^{1,2}$ can be chosen arbitrarily in $[0,1]$.

Moreover, we reformulate the coupling conditions (\ref{eq:merge_kin_eq1}-\ref{eq:merge_kin_eq3}) in terms of $\rho$ and $q$ and obtain $q^3 = q^1+q^2$ and 
 	\begin{align*}
1-(\rho^{1/2}-q^{1/2})&=(1-(\rho^3-q^3))-q^{2/1}
\end{align*}	
or
\begin{align*}
\rho^{1/2}+q^3&=\rho^3+q^{2/1}+q^{1/2}=\rho^3+q^{1}+q^{2}
\end{align*}
and then
\begin{align*}
\rho^{1}&=\rho^2=\rho^{3}\\	
q^3 &= q^1+q^2\ .
\end{align*}		
This is the balance of fluxes and the equality of density.

\subsection{Merging with a priority lane}
\label{sec:prio}
Now we consider a different type of junction. 
In this case lane $1$ has priority and cars on lane $2$ have to wait entering the junction. 
First we consider a situation where lanes $1$ and $3$ run at full capacity.
All cars in road $1$ are driving or $f_0^1=0$. 
Moreover, lane $3$ is completely filled with driving and non driving cars
$f_0^3+f_1^3 = 1$.
In terms of the characteristic quantities this is
\begin{align*}
1-w^1 =1
\end{align*}
and the second condition gives $\frac{q^3}{1-w^3} =1 $ or 
\begin{align*}
Z^3 =1\ .
\end{align*}
Note that the fluxes on road 1 and 3 are given by $q^1=\hat Z^1 (1-w^1)$ and $q^3= Z^3 (1- \hat w^3)$.
Thus the above choice of $w^1$ and $Z^3$ gives the maximal flux given $\hat Z^1 $ and $\hat w^3$.

Mass conservation gives the third condition 
\begin{align*}
f_1^3 = q^3=  q^1+q^2= f_1^1+f_1^2\ .
\end{align*}
Considering this equation in characteristic variables  $w$ and $Z$ and using the above equations, we have 
\begin{align*}
\hat Z^1 + \hat Z^2 (1-w^2) = 1-\hat w^3\ .
\end{align*}

This is directly solved and gives
\begin{align*}
1-w^{2}&= \frac{1-\hat w^3 -  \hat Z^{1}}{\hat Z^{2}}\ .
\end{align*}	

These conditions are well defined,  if $0 \le 1-w^2 \le 1$ and $\hat Z^2 \neq 0$.
This gives the two conditions
\begin{align*}
 \hat w^3 +  \hat Z^{1} \le 1\qquad \text{and}\qquad 1- \hat Z^{2}  \le \hat w^3 +  \hat Z^{1} \ .
\end{align*}	
These conditions have to be truncated to stay within the physical bounds. 
Considering the boundaries of the above domain, see Figure \ref{fig:prio}, we have for $\hat w^3 +  \hat Z^{1} = 1$ that $1-w^2 =0$.
For $\hat w^3 +  \hat Z^{1} = 1-\hat Z^{2}$ we obtain  $1-w^2 =1$.

In case $ \hat w^3 +  \hat Z^{1} \ge 1$ we choose as coupling conditions

\begin{align*}
1- w^2 =0, \;
Z^3=1, \;
1-w^1= \frac{1-\hat w^3 }{ \hat Z^{1} }\ .
\end{align*}	
The last equation is the balance of fluxes with the constraints $1- w^2 =0,
Z^3=1$.

In case $1- \hat Z^{2}  \ge \hat w^3 +  \hat Z^{1} $ we choose as coupling conditions

\begin{align*}
1- w^1 =1, \;
1-w^2= 1, \;
Z^3 = \frac{\hat Z^1+\hat Z^2}{1-\hat w ^3}\,
\end{align*} 
where the last equation is again the balance of fluxes.
\subsubsection{Summary}	

\begin{figure}[h!]
	\begin{center}
		\externaltikz{Statespace_H1}{
			\begin{tikzpicture}[scale = 4]
			\def\dr{0.2}
			\node[below] at (0,0) {$(0,0)$};
			\node[below] at (1,0) {$(1,0)$};
			\node[right] at (1,1) {$(1,1)$};
			\node[left] at (0,1) {$(0,1)$};
			\node[below] at (0.2,0.8) {$\bf I$};
			\node[below] at (0.8,0.5) {$\bf II$};
			\node[right] at (0.03,0.12) {$\bf III$};
			\node[below] at (0.6,0) {$1-\hat w^3$};
			\node[left] at (0,0.6) {$1-\hat w^3$};
			\draw[line width=1pt](0.6,0)--(0,0.6);
			\draw[->](0,0)--(1.2,0) node[below]{$\hat Z^1$};
			\draw[->](0,0)--(0,1.2) node[left]{$\hat Z^2$};
			\draw(1,0)--(1,1);
			\draw(0,1)--(1,1);
			\draw[line width=1pt](0.6,0)--(0.6,1);
			\end{tikzpicture}
		}
	\end{center}
	\caption{Domain of definition for coupling conditions for situation with priority lane.}
	\label{fig:prio}
\end{figure}
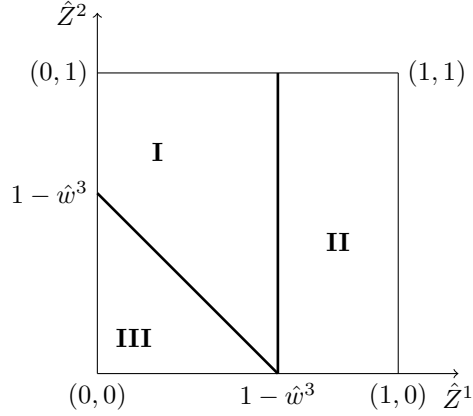

In summary  we have the following  coupling conditions in characteristic variables distinguishing three cases, see Figure \ref{fig:prio}.

{\bf Case I:} $1- \hat Z^{2}  \le  \hat w^3 +  \hat Z^{1} \le 1$ 
\begin{align*}
1-w^1 =1, \;
Z^3 =1, \;
1-w^{2}= \frac{1-\hat w^3 -  \hat Z^{1}}{\hat Z^{2}}\ .
\end{align*}	

{\bf Case II:} $ \hat w^3 +  \hat Z^{1} \ge 1$
\begin{align*}
1- w^2 =0, \;
Z^3=1, \;
1-w^1= \frac{1-\hat w^3 }{ \hat Z^{1} }\ .
\end{align*}

{\bf Case III:} $1- \hat Z^{2}  \ge  \hat w^3 +  \hat Z^{1} $
\begin{align*}
1- w^1 =1, \;
1-w^2= 1, \;
Z^3 = \frac{\hat Z^1+\hat Z^2}{1-\hat w ^3}\ .
\end{align*}

One observes that the restrictions in the above cases in characteristic variables do lead to straightforward restrictions when writing them in macroscopic variables.
Additionally,  we have the physical constraints $0 \le q^1 \le \rho^1\le 1$.
Writing the inequalities in macroscopic variables and using the above conditions in the respective cases, one obtains 

{\bf Case I:} 
$$
\rho^1 = q^1 \le \rho^3 =1,\ q^2 \le \rho^2,\ q^1=q^2+q^3\ ,
$$

{\bf Case II:}
$$
q^2=0 \le q^1 =q^3 \le  \rho^1 \le  \rho^3 = \rho^2=1\ ,
$$

{\bf Case III:}  	
$$
q^1 =\rho^1,\ \rho^2 =q^2  \le q^3 \le \rho^3\ .
$$

Before deriving from these kinetic coupling conditions coupling conditions for the macroscopic equations, we reconsider the kinetic derivation
of boundary conditions for the limit conservation law, compare \cite{BK18} for the present model.

\begin{remark}[Truncated conditions]
\label{sec:merge_priority_delta}
Alternatively the following  coupling conditions  with 
  $0 < \delta <  \bar \delta < 1$ can be considered.

{\bf Case I:} $ (1-\hat w^3)(1-\delta) \ge   \hat Z^{1}$ and $\hat Z^{1} + \hat Z^{2}  \ge  (1-\hat w^3)(1-\delta)  $
\begin{align*}
1-w^1 =1, \;
Z^3 =1-\delta\\
1-w^{2}= \frac{(1-\hat w^3)(1-\delta) -  \hat Z^{1}}{\hat Z^{2}}.
\end{align*}	

{\bf Case II:} $ (1-\hat w^3)(1-\delta) \le   \hat Z^{1}$
\begin{align*}
1- w^2 =0, \;
Z^3=1-\delta\\
1-w^1= \frac{(1-\hat w^3)(1-\delta) }{ \hat Z^{1} }
\end{align*}

{\bf Case III:} $\hat Z^{1} + \hat Z^{2}  \le  (1-\hat w^3)(1-\delta)  $
\begin{align*}
1- w^1 =1, \; 
1-w^2= 1\\
Z^3 = \frac{\hat Z^1+\hat Z^2}{1-\hat w ^3}.
\end{align*}
The capacity in road 3 is large enough to yield  the same macroscopic conditions as the original non-truncated  conditions,
if $Z^3 = 1-\delta $ is chosen such that it is greater than the maximal macroscopic flux, which is given by the maximum of
$ Z(\rho) = \frac{F(\rho)}{1-\rho + F(\rho)} $. This results in 
 $(1-\delta)  \ge Z(\rho) = \frac{F(\rho)}{1-\rho + F(\rho)} $ for all values of $\rho$.
Using the properties of $F$ we observe easily that  we have to choose  $\delta$  as $0\le \delta \le  \frac{1}{1-F^\prime (1) } =\bar \delta$.
\end{remark}

\section{Boundary conditions for the macroscopic equations derived from the kinetic equation}
\label{boundaryconditions}
In this section we determine boundary conditions for the scalar conservation law  from the  boundary value problem of the nonlinear kinetic relaxation system and the kinetic boundary conditions given at the end of section \ref{equations}.
The boundary conditions for the limit equation  are obtained from the kinetic boundary conditions
considering a  half-space problem at the boundary. We refer to  \cite{BSS84,BLP79,CGS,N99} for boundary layers of kinetic equations  and to \cite{AM04,LX96,NT01,WX99,WY99} for investigations of boundary layers for hyperbolic relaxation systems and kinetic equations.

The general procedure is as follows: a half space problem is determined by a rescaling $x\rightarrow \frac{x}{\epsilon}$ of the spatial coordinate in the boundary layer. 
The boundary condition for the half-space problem is given by the original kinetic boundary condition. 
The boundary condition for the limit equation is found by a combination of the asymptotic value of the half-space problem at infinity and a half Riemann problem in the macroscopic variables.

In the following we investigate first the kinetic layer equations and their asymptotic states and then use these results to determine the boundary conditions for the macroscopic problem.

\subsection{Layer solutions for the kinetic  equations}
\label{kinlayer}

The following investigation of the kinetic layer will also be used in the derivation  of the macroscopic coupling conditions
from the kinetic ones.
\subsubsection{Left layer}
Let the left boundary of the domain under consideration be located at $x=x_L$.
Starting from equation \eqref{macro0} and rescaling space as $y= \frac{x-x_L}{\epsilon}$ and neglecting higher order terms in $\epsilon$ one obtains 
the kinetic layer equations  for the left boundary  for $(\rho_L,q_L)$ and $y \in [0, \infty)$ as 
\begin{align}
\label{layerproblem}
\begin{aligned}
\partial_y q_L &=0\\
\frac{q_L}{1-\rho_L}  \partial_y \rho_L  + (1-\frac{q_L}{1-\rho_L})\partial_y  q_L & =- \left(q_l-F(\rho_L) \right) \ .
\end{aligned}
\end{align}
This yields 
\begin{align*}
\begin{aligned}
q_L &=C\\
\partial_y \rho_L   &= (1-\rho_L) \frac{F(\rho_L)-C}{C}\ .
\end{aligned}
\end{align*}
For $0<C < F(\rho^{\star})=\sigma$, where $\rho^{\star}$ denotes the point where the maximum of $F$ is attained,
the above problem has two relevant  fix-points
$$\rho_{-} (C)\le \rho^{\star}\quad  \text{and}\quad  
\rho_+ (C)= \tau (\rho_-) \ge \rho^{\star}\ .$$
Here, $\tau(\rho)\neq \rho$ is defined by $F(\tau(\rho))= F(\rho)$.
$\rho_-$ is instable, $\rho_+ $ is stable. The domain of attraction of the stable fixpoint $\rho_+$  is   the interval $(\rho_-,1)$.

The third fixpoint   $\rho=1$ is not relevant for the further matching procedure, since  it requires $C=0$ in the macroscopic limit. In case  $C=0$ we have the instable fixpoint $\rho_+ = 1$ and the stable fixpoint $\rho_- =0$ with domain of attraction $[0,1)$. 
For $C=F(\rho^{\star})$ we have $\rho_- = \rho_+ = \rho^{\star}$ and 
all solutions with initial values above $\rho^{\star}$ converge towards $\rho^{\star}$, all other solutions  diverge.

\begin{remark}
	In case of the LWR model with $F(\rho) = \rho(1-\rho)$ we have $$\rho_{\pm} (C) = \frac{1}{2} (1 \pm \sqrt{1-4 C})\ ,$$
	with $C< \frac{1}{4}$.
	For $C=\frac{1}{4}$ we have $\rho_- = \rho_+ = \frac{1}{2}$. Moreover, $\tau(\rho) = 1-\rho$, as shown in Figure \ref{figfund}.
	\begin{figure}[h]
		\center
		\externaltikz{LWR}{
				\begin{tikzpicture}[scale = 3]
				\def \rhobar {0.3}	
				\def \rhostar {0.5}
				\draw[->] (0,0)--(1.2,0) node[below]{$\rho$};
				\draw[->] (0,0)--(0,1.2) node[left]{$F(\rho)$}node at (0.2,0.97) {$F(\rho_-)$};
				\draw[black,line width=1pt,domain=0.0:1,smooth,variable=\x,] plot ({\x},{4*\x*(1-\x)}) ;
				\draw[dashed] (\rhobar,{4*\rhobar*(1-\rhobar)})--(\rhobar,0) node[below]{$\rho_-$};
				\draw[dashed] (\rhostar,{4*\rhostar*(1-\rhostar)})--(\rhostar,0) node[below]{$\rho^*$}node at (0.5,1.1) {$\sigma$};;
				\draw[dashed] (\rhobar,{4*\rhobar*(1-\rhobar)})--(1,0)
				node[below]{$1$};	;
				\end{tikzpicture}
				\hspace{0.5cm}
			\begin{tikzpicture}[scale = 3]	
			\draw[->] (0,0)--(1.2,0) node[below]{$C$};
			\draw[->] (0,0)--(0,1.2) node[left]{$\rho$};
			\draw[dashed] (0,0.5)--(1.2,0.5) ;
			\draw[dashed] (1.0,0)--(1.0,1.2) ;
			 \node at (0.5,0.25) {$\rho_-$};
			  \node at (0.5,0.75) {$\rho_+$};
			 \node[left] at (0,0.5) {$\rho^\star$}; 
			 \node[below] at (1.0,0) {$\sigma$}; 
			\draw[black,line width=1pt,domain=0.0:1.0,smooth,variable=\x,samples = 100] plot ({\x},{0.5*(1-sqrt(1-\x))});
			\draw[black,line width=1pt,domain=0.0:1.0,smooth,variable=\x,samples = 100] plot ({\x},{0.5*(1+sqrt(1-\x))});
			\end{tikzpicture}
		}
		\caption{Fundamental diagram, $F(\rho)$ and $\rho_\mp$.}
		\label{figfund}
	\end{figure}
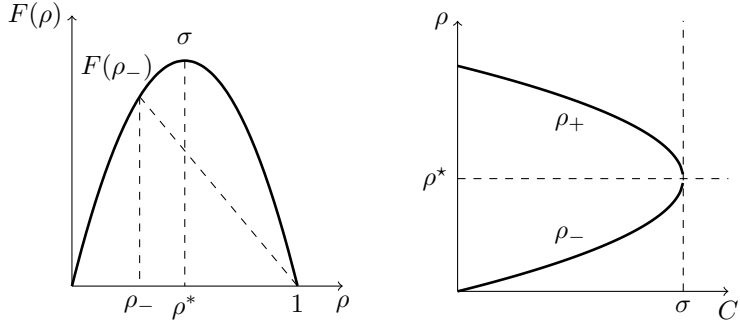
\end{remark}

\subsubsection{Right  layer}

For the right boundary at $x_R$  a scaling $y=\frac{x_R-x}{\epsilon}$ gives the layer equations for $(\rho_R,q_R)$ and  $y\in[0, \infty)$ as
\begin{align}\label{layerright}
\begin{aligned}
q_R &=C\\
- \partial_y \rho_R   &= (1-\rho_R) \frac{F(\rho_R)-C}{C}\ .
\end{aligned}
\end{align}
For  $0 < C < F(\rho^{\star})$
the above problem has again two  relevant fix points
$$\rho_{-}(C) \le \rho^{\star}\ ,\ 
\rho_+ (C) = \tau (\rho_-) \ge \rho^{\star}\; .$$
In this case 
$\rho_- $ is stable, $\rho_+ $ is instable. 
The domain of attraction of the  stable fixpoint $\rho_-$ is $[0,\rho_+)$.

For $C=F(\rho^{\star})=\sigma$ we have $\rho_-= \rho_+ = \rho^{\star}$ and 
all solutions with initial values below  $\rho^{\star}$ converge towards $\rho^{\star}$, all other solutions converge to not admissible states.
For $C=0$ we have the instable fixpoint $\rho_+ = 1$ and the stable fixpoint $\rho_- =0$ with domain of attraction $[0,1)$.

\subsubsection{Summary}
\label{summary}
In summary we have the following cases denoting with  $U$ the unstable fixpoints and with $S$ the stable ones.
Moreover, we use the notation $\rho_K$ for the  values $\rho^\infty_L$ and $\rho^\infty_R$ at infinity of the respective  layers 
and the notation $\rho_0$ for the respective values at $y=0$, i.e. $\rho_L(0)$ and $\rho_R(0)$.
\\
\paragraph{Layer Problem at the left boundary} 
\begin{align*}
&\left.\begin{array}{lll}
\rho_K = \rho_-(C) \quad &\Rightarrow\quad \rho_0 =\rho_-(C), &0 \le   C< \sigma
\end{array}\right\}
&\quad \text{(U)}\\
&\left.\begin{array}{lll}
\rho_K = \rho_+(C) \quad &\Rightarrow\quad \rho_0 \in (\rho_-(C),1), &0< C < \sigma\\
\rho_K = \rho^\star \quad &\Rightarrow\quad  \rho_0 \in [\rho^\star,1),& C=\sigma\\
\rho_K = 1 \quad &\Rightarrow\quad  \rho_0 \in (0,1],& C=0
\end{array}\right\}
&\quad \text{(S)}
\end{align*}

\paragraph{The Layer Problem at the right boundary} 			
\begin{align*}
&\left.\begin{array}{lll}
\rho_K = \rho_+(C) \quad &\Rightarrow\quad \rho_0 =\rho_+(C), &0 \le  C< \sigma 
\end{array}\right\}
&\quad \text{(U)}\\
&\left.\begin{array}{lll}
\rho_K = \rho_-(C) \quad &\Rightarrow\quad \rho_0 \in [0,\rho_+(C)), &0 < C < \sigma\\
\rho_K = \rho^\star \quad &\Rightarrow\quad  \rho_0 \in [0,\rho^\star],& C=\sigma\\
\rho_K = 0 \quad &\Rightarrow\quad  \rho_0 \in [0,1),& C=0
\end{array}\right\}
&\quad \text{(S)}
\end{align*}

We use for the three cases  of the stable fixpoint (S) the notation
$$
\rho_K = \rho_+(C) \quad \Rightarrow\quad  \rho(0) \in \lceil \rho_-(C),1 \rfloor, 0 \le  C \le  \sigma
$$
for the left boundary and 
$$
\rho_K = \rho_-(C) \quad \Rightarrow\quad  \rho(0) \in \lceil 0, \rho_+(C) \rfloor, 0 \le  C \le  \sigma
$$
for the right boundary.

\subsection{Half-Riemann problems for the limit conservation law}
\label{Riemann}
We consider the limit conservation law  $ \partial_t \rho + \partial_x F(\rho)=0$ fulfilling the conditions above.
The possible states $\rho_K$ and $\rho_B $ of the Half-Riemann problems with ingoing waves (shocks and rarefaction waves) at left and right boundary are summarized in the following:\\

\paragraph{The half-Riemann Problem at the left boundary}

\begin{align*}
\rho_B&\leq \rho^\star \ (\text{RP 1})  \quad &\Rightarrow\quad \rho_K &\in [0,\rho^\star ]\\
\rho_B&> \rho^\star \ (\text{RP 2})  \quad &\Rightarrow\quad \rho_K &\in [0,\tau(\rho_B)]\cup\{\rho_B\}
\end{align*}		

\paragraph{The half-Riemann Problem at the right boundary}		
\begin{align*}
\rho_B&\geq  \rho^\star \ (\text{RP 1})  \quad &\Rightarrow\quad \rho_K &\in [ \rho^\star ,1]\\
\rho_B&<  \rho^\star  \ (\text{RP 2})  \quad &\Rightarrow\quad \rho_K &\in \{\rho_B\}\cup [\tau(\rho_B),1]
\end{align*}

The consideration in the last two subsection will now be combined to determine the boundary value as well as coupling conditions for macroscopic equations from the corresponding kinetic conditions.

\subsection{Boundary conditions for the scalar conservation law}
The boundary conditions for the scalar problem are now derived from the kinetic ones by considering the layer equations in section \ref{kinlayer} and the solutions of the half-Riemann problems in section \ref{Riemann}. 
The boundary conditions for the kinetic layer equations (\ref{layerproblem}) at $x=x_L$ are given by the kinetic conditions from section \ref{equations}. 
In the layer equations for $(\rho_L,q_L)$ at the left boundary we prescribe $Z$ at $y=0$, for the layer equations at the right boundary for $(\rho_R,q_R)$ we prescribe $w$ at $y=0$. 
Then one determines the asymptotic states  $\rho^\infty_L$ and $\rho^\infty_R$ at $y=\infty$, i.e. the solution at infinity of the layer equations.
This state is finally used as boundary condition for the scalar equations. 
The initial trace at the boundary of the scalar equation is in the following denoted by $\rho_B$.

\subsubsection{Left boundary}
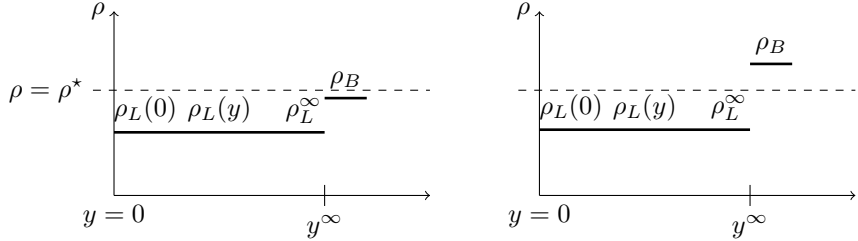
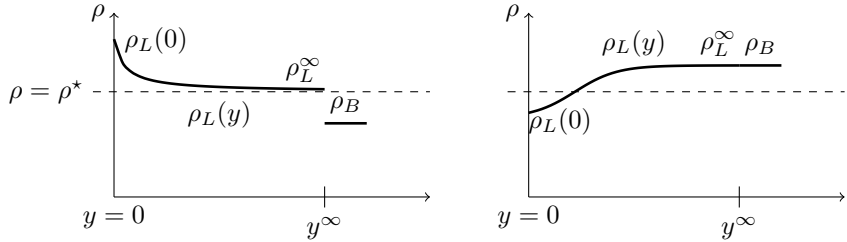
\begin{figure}[h]
	\begin{subfigure}[b]{0.49\textwidth}
		\externaltikz{Layer_cases1}{
			\begin{tikzpicture}[scale = .7]
			\def\len{4}
			\def\low{-2}
			
			\def \rhozero {-1}
			\def \rhoB {-0.15}
			\def \rhoK {-0.8}
			\def \rhok {-0.75}
			\def \Ckonst {\rhok*\rhok}
			\def \Ctwo {-1/(2*\rhok)*ln((\rhok+\rhozero)/(\rhok-\rhozero))}
			
			\node[] (xinf) at (\len,\low){};
			\node[below] at (0,\low){$y=0$};
			\draw[dashed] (-0.1*\len,0)--(1.5*\len,0) node[left,pos = 0]{$\rho=\rho^\star$};
			\draw (0,\low)--(\len,\low);
			\draw[->] (\len,\low)--(1.5*\len,\low);
			\draw[->] (0,\low)--(0,1.5) node[left]{$\rho$};
			\draw (\len,0.9*\low)--(\len,1.1*\low) node[below]{$y^\infty$};
			\draw[line width=1pt] (\len,\rhoB)--(1.2*\len,\rhoB) node[above,pos=0.5] {$\rho_B$};
			\draw[line width=1pt,domain=0.0:\len,smooth,variable=\x,] plot ({\x},{\rhoK}) ;
			\node[above] at (0.15*\len,\rhoK) {$\rho_L(0)$};
			\node[above] at (0.5*\len,\rhoK) {$\rho_L(y)$};
			\node[above] at (\len*0.9,\rhoK) {$\rho^\infty_L$};
			\end{tikzpicture}
			
		}
		\caption{Layer and RP solution for case 1 a). }
	\end{subfigure}
	\begin{subfigure}[b]{0.49\textwidth}
		\externaltikz{Layer_cases2}{	
			\begin{tikzpicture}[scale = .7]
			\def\len{4}
			\def\low{-2}	
			\def \rhozero {-1}
			\def \rhoB {0.5}
			\def \rhok {-0.75}
			\def \Ckonst {\rhok*\rhok}
			\def \Ctwo {-1/(2*\rhok)*ln((\rhok+\rhozero)/(\rhok-\rhozero))}
			\node[] (xinf) at (\len,\low){};
			\node[below] at (0,\low){$y=0$};
			\draw[dashed] (-0.1*\len,0)--(1.5*\len,0);
			\draw (0,\low)--(\len,\low);
			\draw[->] (\len,\low)--(1.5*\len,\low);
			\draw[->] (0,\low)--(0,1.5) node[left]{$\rho$};
			\draw (\len,0.9*\low)--(\len,1.1*\low) node[below]{$y^\infty$};
			\draw[line width=1pt] (\len,\rhoB)--(1.2*\len,\rhoB) node[above,pos=0.5] {$\rho_B$};
			\draw[line width=1pt,domain=0.0:\len,smooth,variable=\x,] plot ({\x},{\rhok}) ;
			\node[above] at (0.15*\len,\rhok) {$\rho_L(0)$};
			\node[above] at (0.5*\len,\rhok) {$\rho_L(y)$};
			\node[above] at (\len*0.9,\rhok) {$\rho_L^\infty$};	
			\end{tikzpicture}
		}
		\caption{Layer and RP solution for case 1 b). }
	\end{subfigure}
	\begin{subfigure}[b]{0.49\textwidth}
		\externaltikz{Layer_cases3}{
			\begin{tikzpicture}[scale = .7]
			\def\len{4}
			\def\low{-2}
			
			\def \rhozero {-1}
			\def \rhoB {-0.6}
			\def \rhok {0}
			\def \Ckonst {\rhok*\rhok}
			\def \Ctwo {-1/(2*\rhok)*ln((\rhok+\rhozero)/(\rhok-\rhozero))}
			
			\node[] (xinf) at (\len,\low){};
			\node[below] at (0,\low){$y=0$};
			\draw[dashed] (-0.1*\len,0)--(1.5*\len,0) node[left,pos = 0]{$\rho = \rho^\star$};
			\draw (0,\low)--(\len,\low);
			\draw[->] (\len,\low)--(1.5*\len,\low);
			\draw[->] (0,\low)--(0,1.5) node[left]{$\rho$};
			\draw (\len,0.9*\low)--(\len,1.1*\low) node[below]{$y^\infty$};
			\draw[line width=1pt] (\len,\rhoB)--(1.2*\len,\rhoB) node[above,pos=0.5] {$\rho_B$};
			\draw[line width=1pt,domain=0.0:\len,smooth,variable=\x,] plot ({\x},{1/(\x/\len*20-1/(\rhozero))}) ;
			\node[below] at (0.2*\len,1.4) {$\rho_L(0)$};
			\node[below] at (0.5*\len,\rhok) {$\rho_L(y)$};
			\node[above] at (\len*0.9,\rhok) {$\rho_L^\infty$};	
			\end{tikzpicture}
		}
		\caption{Layer and RP solution for case 2. }
	\end{subfigure}
	\begin{subfigure}[b]{0.49\textwidth}
		\externaltikz{Layer_cases4}{
			\begin{tikzpicture}[scale = .7]
			\def\len{4}
			\def\low{-2}
			
			\def \rhozero {0.4}
			\def \rhoB {0.5}
			\def \rhok {0.5}
			\def \Ckonst {\rhok*\rhok}
			\def \Ctwo {-1/(2*\rhok)*ln((\rhok+\rhozero)/(\rhok-\rhozero))}
			
			\node[] (xinf) at (\len,\low){};
			\node[below] at (0,\low){$y=0$};
			\draw[dashed] (-0.1*\len,0)--(1.5*\len,0);
			\draw (0,\low)--(\len,\low);
			\draw[->] (\len,\low)--(1.5*\len,\low);
			\draw[->] (0,\low)--(0,1.5) node[left]{$\rho$};
			\draw (\len,0.9*\low)--(\len,1.1*\low) node[below]{$y^\infty$};	
			\draw[line width=1pt] (\len,\rhoB)--(1.2*\len,\rhoB) node[above,pos=0.5] {$\rho_B$};
			\draw[line width=1pt,domain=0.0:\len,smooth,variable=\x,] plot ({\x},{-\rhok*tanh(-\rhok*(\x/\len*10+\Ctwo))}) ;
			\node[above] at (0.15*\len,-1.0) {$\rho_L(0)$};
			\node[above] at (0.5*\len,\rhoB) {$\rho_L(y)$};
			\node[above] at (\len*0.9,\rhoB) {$\rho_L^\infty$};	
			\end{tikzpicture}
		}	
		
		\caption{Layer and RP solution for case 3. }
	\end{subfigure}
	\caption{The solution of the boundary layer $\rho_L$ and the jump due to Riemann problem at $y^\infty$ at the left boundary of the domain.}
	\label{layerfig}
\end{figure}

Assume for the left boundary  $0 \le Z(x_L)$ to be known and consider the layer equation (\ref{kinlayer}) for $y \in [0, \infty)$ with the left boundary condition $ \frac{ q_L(0)}{1-\rho_L(0)+q_L(0)}= Z(x_L) $. 
We distinguish the following  cases. 
An illustration of the different situations is given in Figure \ref{layerfig}.

\noindent {\bf Case 1: ingoing flow}
We have two sub-cases

\noindent{\bf Case 1.a:}
$\rho_B \le \rho^{\star}$ (RP 1) and unstable layer solution with $\rho_L(0)=\rho_L^\infty = \rho_K=\rho_-(C) \in [0, \rho^\star]$.

\noindent{\bf Case 1.b:}
  $\rho_B > \rho^{\star}$ (RP 2) and unstable layer solution with $ \rho_L(0) = \rho_K =\rho_-(C)  \in [0,\tau(\rho_B)] $.
  
The  layer solution is in both cases the unstable solution
$$
\rho_L(y) = \rho_-(C) \le \rho^{\star}\ ,
$$
where $0 \le \rho_- \le \rho^\star$  is determined from $Z(x_L)
 =\frac{C}{1-\rho_L(0) +C}=  \frac{F(\rho_-)}{1-\rho_- +F(\rho_-)} $.
 In Case 1.a this equation has a unique solution $\rho_-$  if  $Z(x_L)$
fulfills  $0 \le Z(x_L)\le \frac{\sigma}{1-\rho^{\star}+ \sigma}$.
In Case 1.b we have  $\rho_L(0) = \rho_- \le\tau(\rho_B)$, if $  0 \le Z(x_L)\le \frac{F(\tau(\rho_B))}{1-\tau(\rho_B)+F(\tau(\rho_B))} $.
In both cases one obtains for the solution $\rho$ of the scalar problem a boundary condition $\rho(x_L) = \rho_L(0)$ which gives a wave with positive speed starting at the boundary.

\noindent{\bf Case 2: transonic flow} $\rho_B <\rho^{\star}$ (RP 1)  and stable layer solution with 
$\rho_K =\rho_+(C)  \in [0, \rho^\star]$ and $ \rho_L(0) \in (\rho_-(C),1)$.

One has to choose $C$ as the maximal possible value $C = F(\rho^{\star})$.
From $ Z(x_L)= \frac{C}{1-\rho_L(0)+C}= \frac{\sigma}{1-\rho_L(0)+\sigma}$, we obtain $\rho_L(0)$ with
 $\rho_L(0)\ge \rho^{\star}$, if $Z(x_L) \ge  \frac{\sigma}{1-\rho^{\star}+\sigma}$.
The layer solution is no longer constant in space.
Moreover,  $\rho^\infty_L=\rho^{\star}$.
In this case one obtains a rarefaction wave.

\noindent{\bf Case 3: outgoing flow}
$\rho_B >\rho^{\star}$ (RP 2) and stable layer solution with $\rho_K =\rho_+(C) = \rho_B $ and 
$ \rho(0) \in (\rho_-(C),1)$

Here,
$
\rho^\infty_L = \rho_+(C) = \rho_B
$
yields $C$ and 
$Z(x_L)= \frac{C}{1-\rho_L(0)+C}= \frac{F(\rho_+)}{1-\rho_L(0)+F(\rho_+)}= \frac{ F(\rho_B)}{1-\rho_L(0)+F(\rho_B)}$ 
gives
$
\rho_L(0) \ge \tau(\rho_B)\ ,
$
if $Z(x_L)\ge \frac{F(\tau(\rho_B))}{1-\tau(\rho_B)+F(\tau(\rho_B))}$.
There is no wave starting at the boundary and we do not prescribe a boundary condition for the scalar conservation law. The layer does not have a constant solution.

\subsubsection{Right boundary}

For the right boundary we prescribe the 1-Riemann invariant $w(x_R)$ and consider equation (\ref{layerright}) for $y \in [0, \infty)$ with
the boundary condition $  \rho_R(0) - q_R(0)= w(x_R) $. We distinguish again  three cases.

\noindent{\bf Case 1: ingoing flow}

\noindent{\bf Case 1.a:}
$\rho_B \ge \rho^{\star}$ (RP 1) and unstable layer solution with $\rho_R(0)=\rho_K =\rho_+(C) \in [\rho^\star,1]$.

\noindent{\bf Case 1.b:}
$\rho_B < \rho^{\star}$ (RP 2) and unstable layer solution with $ \rho_R(0) = \rho_K =\rho_+(C)  \in [\tau(\rho_B),1] $.

The unstable layer solution is 
$$
\rho_R(y) = \rho_+ (C) \ge \rho^{\star}\ .
$$
Here,  $0<C< F(\rho^{\star})$ is determined from $w(x_R) 
=  \rho_+(C) -C $. We determine $\rho_+$ from 
$\rho_+-F(\rho_+)= w(x_R)$. In the first case this has a unique solution $\rho^+\in [\rho^*,1]$ if $ w(x_R) \ge \rho^{\star}- F(\rho^{\star})$.
In the second case it has a unique solution $\rho^+\in[\tau(\rho_B),1]$ if $w(x_R)\ge \tau(\rho_B)- F(\tau(\rho_B)) $.
Moreover, 
$\rho_R(0) = \rho_+ \ge \rho^{\star}$
and
$\rho^\infty_R=\rho_R(0)\ . $
In the first case $\rho_R(0) \ge \rho^{\star}$ and in the second case $\rho_R(0) \ge\tau(\rho_B)$.

\noindent{\bf Case 2: transonic flow} 
$\rho_B >\rho^{\star}$ (RP 1)  and stable layer solution with 
$\rho_K =\rho_-(C)  \in [ \rho^\star,1]$ and $ \rho(0) \in [0, \rho_+(C))$.
In this case we have $C = F(\rho^{\star})$.
From $w(x_R) = \rho_R(0)-C= \rho_R(0)-F(\rho^{\star})$ we obtain $\rho_R(0) = w(x_R)+ F(\rho^{\star})$.
This yields $\rho_R(0)\le \rho^{\star}$, if  $w(x_R) \le \rho^{\star}-  F(\rho^{\star})$.
Moreover,  $\rho^\infty_R=\rho^{\star}$.

\noindent{\bf Case 3: outgoing flow}
$\rho_B <\rho^{\star}$ (RP 2) and stable layer solution with $\rho_K =\rho_-(C) =\rho_B $ and 
$ \rho(0) \in [0,\rho_+(C))$.
Then,
$
\rho^\infty_R  = \rho_-(C) = \rho_B\ .
$
This yields $C$
and 
$w(x_R)= \rho_R(0)-C= \rho_R(0)- F(\rho_1)= \rho_R(0)- F(\rho_B)$ gives
$$
\rho_R(0) = w(x_R)+ F(\rho_B)\le \tau(\rho_B)\ ,
$$
if $w(x_R)\le \tau(\rho_B)- F(\tau(\rho_B)) $.

\section{Macroscopic coupling conditions for junctions with  fair merging}
\label{macroscopiccc}

We consider the fair merging case.
In order to determine the coupling conditions for the macroscopic equations we proceed analogously as for the boundary conditions in the last section.
First we investigate  the kinetic layers at the nodes coupled to each other via the coupling conditions and determine their asymptotic states. Then, we   match these  results to Riemann solutions of the  macroscopic problems on each of the roads.

Assuming the boundary traces $\rho_B^1,\rho_B^2,\rho_B^3$ on the three roads to be given, we have to determine the new states $\rho_K^1,\rho_K^2$ and $\rho_K^3$ at the node. On the one hand $\rho_K^1,\rho_K^2$ and $\rho_K^3$ are the asymptotic states of the respective layer problems, on the other hand they are the right (for road 1 and 2) or left (for road 3) states of the half-Riemann problems with $\rho_B^1,\rho_B^2,\rho_B^3$
as the corresponding left (road 1 and 2) or right  state (road 3).
The states at the junction (corresponding to $y=0$ for the layers) are denoted in the following and in section \ref{technical} by
$\rho_0$.
We have to consider eight different configurations of Riemann problems.
For  each of them all possible combinations with stable or unstable layer solutions have to be discussed.
Not admissible combinations are not listed.
The proof of the following statement is given  in Section \ref{technical}.
For a detailed discussion of the coupling of the layer solutions we refer to \ref{layerproof} and for the matching of the layer solutions to the half Riemann problems   to
\ref{proof}.

\begin{theorem}[Macroscopic coupling conditions: Fair merging]
\label{theoremmerge}
Starting from the kinetic coupling conditions in section \ref{sec:fair} the asymptotic derivation of the  coupling conditions  for the macroscopic equations gives  the following  cases using the notation 
RP1/2-1/2-1/2 for the respective combination of the half Riemann problems.

\noindent {\bf Case 1, RP1-1-1.} $\rho_B^1 \ge \rho^\star , \rho_B^2 \ge \rho^\star , \rho_B^3 \le  \rho^\star $.
Then $C_1=C_2=\frac{\sigma}{2}, C_3=\sigma$ and
\begin{align*}
\rho_K^1&= \rho_+(\frac{\sigma}{2}) &\rho_K^2&= \rho_+(\frac{\sigma}{2}) &\rho_K^3 =\rho^\star \ .
\end{align*}

\noindent {\bf Case 2, RP1-1-2.}  $\rho_B^1 \ge \rho^\star , \rho_B^2 \ge \rho^\star , \rho_B^3 \ge  \rho^\star $
$C_1=C_2=F(\frac{\rho_B^3}{2}), C_3=F(\rho_B^3)$ and	
\begin{align*}
\rho_K^1&= \rho_+(\frac{F(\rho_B^3)}{2}) &\rho_K^2&= \rho_+(\frac{F(\rho_B^3)}{2}) &\rho_K^3 =\rho_B^3 \ .
\end{align*}

\noindent {\bf Case 3, RP1-2-1} $\rho_B^1 \ge \rho^\star , \rho_B^2 \le \rho^\star , \rho_B^3 \le  \rho^\star $. Then
we have  two cases:

\noindent 
If $F(\rho_B^2) \ge \frac{\sigma}{2}$, then
$C_1=C_2=\frac{\sigma}{2}, C_3=F(\rho_B^3)$ and
\begin{align*}
\rho_K^1&= \rho_+(\frac{\sigma}{2}) &\rho_K^2&= \rho_+(\frac{\sigma}{2}) &\rho_K^3 =\rho^\star \ .
\end{align*}
If $F(\rho_B^2) \le \frac{\sigma}{2}$, then $C_1=\sigma-F(\rho_B^2), C_2=F(\rho_B^2), C_3=\sigma$ and 
\begin{align*}
\rho_K^1&= \rho_+(\sigma-F(\rho_B^2)) &\rho_K^2&= \rho_B^2&\rho_K^3 =\rho^\star \ .
\end{align*}

\noindent {\bf Case 4,	RP2-1-1} $\rho_B^1 \le \rho^\star , \rho_B^2 \ge \rho^\star , \rho_B^3 \le  \rho^\star $. The
case is symmetric to Case 3.

\noindent
If $F(\rho_B^1) \ge \frac{\sigma}{2}$, then $C_1=C_2=\frac{\sigma}{2}, C_3=\sigma$ and
\begin{align*}
\rho_K^1&= \rho_+(\frac{\sigma}{2}) &\rho_K^2&=  \rho_+(\frac{\sigma}{2})&\rho_K^3 =\rho^\star \ .
\end{align*}
If $F(\rho_B^1) \le \frac{\sigma}{2}$, then $C_1=F(\rho_B^1), C_2=\sigma-F(\rho_B^1), C_3=\sigma$ and
\begin{align*}
\rho_K^1&= \rho_B^1&\rho_K^2&= \rho_+(\sigma-F(\rho_B^1)&\rho_K^3 =\rho^\star \ .
\end{align*}

\noindent {\bf Case 5, RP1-2-2} $\rho_B^1 \ge \rho^\star , \rho_B^2 \le \rho^\star , \rho_B^3 \ge  \rho^\star .$ Then 
we have 2 cases.

\noindent
If $F(\rho_B^3) \le 2 F(\rho_B^2)$, then $C_1=\frac{F(\rho_B^3)}{2}=C_2, C_3=F(\rho_B^3)$ and
\begin{align*}
\rho_K^1&= \rho_+(\frac{F(\rho_B^3)}{2}) &\rho_K^2&= \rho_+(\frac{F(\rho_B^3)}{2}) &\rho_K^3 =\rho_B^3\ .
\end{align*}
If $F(\rho_B^3) \ge 2 F(\rho_B^2)$, then $C_1=F(\rho_B^3)- F(\rho_B^2), C^2 = F(\rho_B^2), C_3=F(\rho_B^3)$ and
\begin{align*}
\rho_K^1&= \rho_+(F(\rho_B^3)-F(\rho_B^2)) &\rho_K^2&=  \rho_B^2&\rho_K^3 =\rho_B^3 \ .
\end{align*}

\noindent{\bf Case 6, RP2-1-2}
$\rho_B^1 \le \rho^\star , \rho_B^2 \ge \rho^\star , \rho_B^3 \ge  \rho^\star $.  The case is symmetric to Case 5.

\noindent
If $F(\rho_B^3) \le 2 F(\rho_B^1)$, then $C_1=\frac{F(\rho_B^3)}{2}=C_2, C_3=F(\rho_B^3)$ and
\begin{align*}
\rho_K^1&= \rho_+(\frac{F(\rho_B^3)}{2}) &\rho_K^2&= \rho_+(\frac{F(\rho_B^3)}{2}) &\rho_K^3 =\rho_B^3 .
\end{align*}
If $F(\rho_B^3) \ge 2 F(\rho_B^1)$, then $C_1=F(\rho_B^1), C^2 = F(\rho_B^3)- F(\rho_B^1), C_3=F(\rho_B^3)$ and
\begin{align*}
\rho_K^1&= \rho_B^1&\rho_K^2&=  \rho_+(F(\rho_B^3)- F(\rho_B^1))&\rho_K^3 =\rho_B^3 \ .
\end{align*}

\noindent {\bf Case 7,	RP2-2-1}  $\rho_B^1 \le \rho^\star , \rho_B^2 \le \rho^\star , \rho_B^3 \le  \rho^\star $. We have 4 subcases.

\noindent
If $F(\rho_B^1) +  F(\rho_B^2) \le \sigma$, then $C_1=F(\rho_B^1), C_2=F(\rho_B^2), C_3=F(\rho_B^1)+F(\rho_B^2)$ and
\begin{align*}
\rho_K^1&= \rho_B^1&\rho_K^2&= \rho_B^2 &\rho_K^3 =\rho_-(F(\rho_B^1+F(\rho_B^2))\ .
\end{align*}
If $F(\rho_B^1) +  F(\rho_B^2) \ge \sigma$ and $F(\rho_B^1) \ge \frac{\sigma}{2}, F(\rho_B^2) \ge  \frac{\sigma}{2}$, then 
$C_1= \frac{\sigma}{2}=C^2, C^3 =\sigma $ and
\begin{align*}
\rho_K^1&= \rho_+( \frac{\sigma}{2}) &\rho_K^2&=  \rho_+( \frac{\sigma}{2})&\rho_K^3 =\rho^\star\ .
\end{align*}
If $F(\rho_B^1) +  F(\rho_B^2) \ge \sigma$ and $F(\rho_B^1) \ge \frac{\sigma}{2}, F(\rho_B^2) \le  \frac{\sigma}{2}$, then 
$C_1=\sigma-F(\rho_B^2), C_2=F(\rho_B^2)$, $ C_3=\sigma$ and
\begin{align*}
\rho_K^1&= \rho_+(\sigma-F(\rho_B^2)) &\rho_K^2&= \rho_B^2 &\rho_K^3 =\rho^\star\ .
\end{align*}
If $F(\rho_B^1) +  F(\rho_B^2) \ge \sigma$ and $F(\rho_B^1) \le \frac{\sigma}{2}, F(\rho_B^2) \ge  \frac{\sigma}{2}$,  then 
$C_1=F(\rho_B^1), C_2=\sigma-F(\rho_B^1)$, $C_3=\sigma$ and 
\begin{align*}
\rho_K^1&= \rho_B^1&\rho_K^2&= \rho_+(\sigma- F(\rho_B^1)) &\rho_K^3 =\rho^\star\ .
\end{align*}

\noindent {\bf Case 8, RP2-2-2}
$\rho_B^1 \le \rho^\star , \rho_B^2 \le \rho^\star , \rho_B^3 \ge  \rho^\star $.
We have 4 subcases. 

\noindent
If $F(\rho_B^3) \le 2   F(\rho_B^1) $ and $F(\rho_B^3) \le 2   F(\rho_B^2) $, then $C_1=\frac{F(\rho_B^3)}{2}=C^2,  C_3=F(\rho_B^3)$ and
\begin{align*}
\rho_K^1&= \rho_+(\frac{F(\rho_B^3)}{2}) &\rho_K^2&= \rho_+(\frac{F(\rho_B^3)}{2})  &\rho_K^3 =\rho_B^3 .
\end{align*}
If $F(\rho_B^3) \ge 2   F(\rho_B^2) $ and $F(\rho_B^1)+F(\rho_B^2) \ge   F(\rho_B^3) $, then $C^2=F(\rho_B^2)$, $C^3=F(\rho_B^3)$, $C_1=F(\rho_B^3)-F(\rho_B^2)$ and
\begin{align*}
\rho_K^1&= \rho_+(F(\rho_B^3)-F(\rho_B^2)) &\rho_K^2&= \rho_B^2 &\rho_K^3 =F(\rho_B^3 .
\end{align*}
If $F(\rho_B^3) \ge 2   F(\rho_B^1) $ and $F(\rho_B^1)+F(\rho_B^2) \ge   F(\rho_B^3) $, then $C_1=F(\rho_B^1)$,  $C^3=F(\rho_B^3)$, $C^2=F(\rho_B^3)-F(\rho_B^1)$ and
\begin{align*}
\rho_K^1&= \rho_B^1 &\rho_K^2&= \rho_+(F(\rho_B^3)-F(\rho_B^1))  &\rho_K^3 =\rho_B^3\ .
\end{align*}
If $F(\rho_B^1)+F(\rho_B^2) \le   F(\rho_B^3) $, then $C_1=F(\rho_B^1),C^2=F(\rho_B^2)$, $C^3=F(\rho_B^1)+F(\rho_B^2)$
and 
\begin{align*}
\rho_K^1&=\rho_B^1 &\rho_K^2&= \rho_B^2  &\rho_K^3 =\rho_-(F(\rho_B^1)+F(\rho_B^2))\ .
\end{align*}

\end{theorem}

For the determination of the kinetic densities at the coupling points $\rho_0^i$ we refer to the proof in section \ref{proof}.
Note that their values are not needed for the macroscopic model.
The above  result can be rewritten in a more convenient way.

\subsection{Supply-Demand formulation of the coupling conditions}

\begin{figure}[h]
		\center
		\externaltikz{supplydemand}{
			\begin{tikzpicture}[scale = 3]
				\def \rhobar {0.3}	
				\def \rhostar {0.5}
				\draw[->] (0,0)--(1.2,0) node[below]{$\rho$};
				\draw[->] (0,0)--(0,1.2) node[left]{$F(\rho)$};
				\draw[dashed] (1,1.2)--(1,0.0) node[below]{$1$};
				\draw[black,line width=1pt,domain=\rhostar:1,smooth,variable=\x,] plot ({\x},{1}) ;
				\draw[black,line width=1pt,domain=0.0:\rhostar,smooth,variable=\x,] plot ({\x},{4*\x*(1-\x)}) ;
				\draw[dashed,black,line width=1pt,domain=\rhostar:1,smooth,variable=\x,] plot ({\x},{4*\x*(1-\x)}) ;
				\draw[dashed] (\rhostar,{4*\rhostar*(1-\rhostar)})--(\rhostar,0) node[below]{$\rho^*$}node at (0.5,1.1) {$\sigma$};
			\end{tikzpicture}
			\hspace{0.5cm}
			\begin{tikzpicture}[scale = 3]
				\def \rhobar {0.3}	
				\def \rhostar {0.5}
				\draw[->] (0,0)--(1.2,0) node[below]{$\rho$};
				\draw[->] (0,0)--(0,1.2) node[left]{$F(\rho)$};
				\draw[dashed] (1,1.2)--(1,0.0) node[below]{$1$};
				\draw[black,line width=1pt,domain=0.0:\rhostar,smooth,variable=\x,] plot ({\x},{1}) ;
				\draw[black,line width=1pt,domain=\rhostar:1,smooth,variable=\x,] plot ({\x},{4*\x*(1-\x)}) ;
				\draw[dashed,black,line width=1pt,domain=0.0:\rhostar,smooth,variable=\x,] plot ({\x},{4*\x*(1-\x)}) ;
				\draw[dashed] (\rhostar,{4*\rhostar*(1-\rhostar)})--(\rhostar,0) node[below]{$\rho^*$}node at (0.5,1.1) {$\sigma$};;
			\end{tikzpicture}
		
		}
	
		\caption{Supply- and demand functions $c^i$ for ingoing (left) and outgoing (right) roads.}
		\label{fig:supply}
	\end{figure}
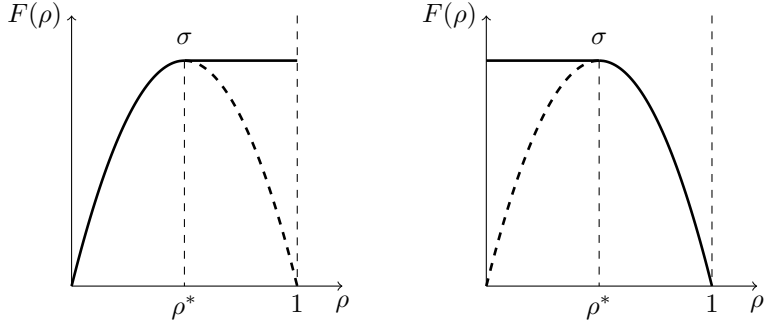

We use the supply-demand representation and denote the sets of valid resulting fluxes $C^i$ by $\Omega^i$, compare \cite{CGP05,L,Dag1,Dag2,HR} and Figure \ref{fig:supply}. 
For the incoming roads $i=1,2$ this is
\begin{align*}
	 \rho_B^i \le \rho^\star \Rightarrow \Omega^i = [0, F(\rho_B^i)] &&\text{ and }&&
	\rho_B^i \ge \rho^\star \Rightarrow \Omega^i = [0, \sigma]\ .	 
\end{align*}
For the outgoing road $i=3$
\begin{align*}
\rho_B^i \le\rho^\star \Rightarrow \Omega^i = [0, \sigma]  &&\text{ and }&&
 \rho_B^i\ge  \rho^\star \Rightarrow \Omega^i = [0, F(\rho_B^i)]\ .	 
\end{align*}
We define the maximal admissible flux $c^i$ such that $\Omega^i = [0, c^i]$.
Using this notation the above 8 cases are rewritten as  

{\bf Case 1, RP1-1-1.} 
This is a case with  $c^1,c^2\le \frac{c^3}{2}: C^1 =c^1, C^2 = c^2.$

{\bf Case 2, RP1-1-2.}  
This is a case with $c^1,c^2\ge \frac{c^3}{2}: C^1 =C^2 = \frac{c^3}{2}.$

{\bf Case 3, RP1-2-1} 
We have $c^1 \ge c^2$ and  two cases:
\begin{align*}
c^2 \ge \frac{c^3}{2}&: C^1= C^2 = \frac{c^3}{2}\\
c^2 \le \frac{c^3}{2}&:  C^1= c^3-c^2, C^2 =c^2.
\end{align*}

{\bf Case 4,	RP2-1-1} Symmetric to Case 3.
We have $c^1 \le c^2$ and  two cases:
\begin{align*}
c^1 \ge \frac{c^3}{2}&: C^1= C^2 = \frac{c^3}{2}\\
c^1 \le \frac{c^3}{2}&:  C^1= c^1, C^2 =c^3 -c^1.
\end{align*}

{\bf Case 5, RP1-2-2} 
In terms of the  $c^i $ this case is the same as Case 3.

{\bf Case 6, RP2-1-2}
This case is the same as Case 4.

{\bf Case 7,	RP2-2-1}  
We have four cases:
\begin{align*}
c^1 +c^2 \le c^3&:C^1= c^1, C^2 = c^2\\
c^1 +c^2 \ge c^3,c^1 \ge \frac{c^3}{2}, c^2\ge \frac{c^3}{2} &:C^1= C^2 =\frac{c^3}{2}\\
c^1 +c^2 \ge c^3,c^1 \ge \frac{c^3}{2}, c^2\le \frac{c^3}{2}&:C^1= c^3-c^2,C^2 =c^2\\
c^1 +c^2 \ge c^3,c^1 \le \frac{c^3}{2}, c^2\ge \frac{c^3}{2}  &:C^1= c^1, C^2 =c^3-c^1.
\end{align*}

{\bf Case 8, RP2-2-2}
We obtain the same as in Case 7.\\

\noindent All in all, only 4 different cases are left:

	{\bf Case A:}
		\begin{align*}
		c^1 +c^2 \le c^3:C^1= c^1, C^2 = c^2\ ,
		\end{align*}
	\qquad {\bf Case B:}
		\begin{align*}
		c^1 +c^2 \ge c^3,c^1 \ge \frac{c^3}{2}, c^2\ge \frac{c^3}{2}:  C^1= C^2 =\frac{c^3}{2}\ ,
		\end{align*}
	\qquad {\bf Case C:}
		\begin{align*}
		c^1 +c^2 \ge c^3,c^1 \ge \frac{c^3}{2}, c^2\le \frac{c^3}{2}:  C^1= c^3-c^2,C^2 =c^2\ ,
		\end{align*}
	\qquad {\bf Case D:}
		\begin{align*}
		c^1 +c^2 \ge c^3,c^1 \le \frac{c^3}{2}, c^2\ge \frac{c^3}{2}:  C^1= c^1, C^2 =c^3-c^1\ .
		\end{align*}

We note that these can be rewritten as
\begin{align}\label{eq:LWRcoupling}
\begin{aligned}
	C^i &= 
	\begin{cases}
	c^i& \text{ if } c^1+c^2\leq c^3\\
	\min\left(c^i,c^3-\min\left(c^1,c^2,\frac{c^3}{2}\right)\right)
	&\text{ else }
	\end{cases}
	\qquad i = 1,2\ ,\\
	C^3 &= C^1+C^2\ .
\end{aligned}
\end{align}

That means we merge symmetric, if both incoming roads have a  flux which is larger than their share  in the outgoing road. Otherwise we reduce only the flux in the road with larger flux.
We refer, for example, to \cite{GPBook} for similar  coupling conditions for scalar conservation laws.

\begin{remark}
	The above derivation shows that a classical merge condition as the above for a nonlinear scalar conservation law
	can be interpreted on the kinetic level as a combination of the balance of fluxes and the equality of densites on all roads.
	We note once more that the equality of densities is not fulfilled on the macroscopic level of the conservation law,
	only the balance of fluxes is common for both levels of coupling conditions.
\end{remark}

\subsection{Macroscopic coupling conditions for situation with priority lane}
\label{sec:macro_coupling_priority}

We use the same  notation as in the previous section, i.e. we
define $C^i,c^i$ and the sets $\Omega^i=[0, c^i]$ as above, depending on
whether incoming or outgoing roads are considered.
The kinetic  conditions for situations wit a priority lane  from Section \ref{sec:prio} lead 
in the limit to the macroscopic coupling conditions

{\bf Case A:} 
\begin{align*}
c^1 +c^2 \le c^3:C^1= c^1, C^2 = c^2\ ,
\end{align*}
\qquad {\bf Case B:}
\begin{align*}
c^1 +c^2 \ge c^3,c^1 \ge c^3:  C^1= c^3, C^2= 0\ ,
\end{align*}
\qquad {\bf Case C:}
\begin{align*}
c^1 +c^2 \ge c^3,c^1 \le c^3:  C^1= c^1,C^2 =c^3-c^1\ ,
\end{align*}
The last two cases  can be written as 
\begin{align*}
 C^1 &= \min\Big(  c^1, c^3 \Big)\\
C^2 &= c^3-C^1= \max \Big(c^3-c^1,0 \Big) \\
C^3&=C^1+C^2\ .
\end{align*}
This means, as long as the capacity of road 3 is large enough, the flow in road 1 is unrestrained.

\begin{remark}
More generally one considers macroscopic merge conditions of the following form 
 \cite{CGP05,HK,L,Dag1,Dag2}:
As in the two cases considered here, one uses  $C^1=c^1$ and $C^2=c^2$ and $C^3 = C^1+ C^2$, if  $c^1 +c^2 \le  c^3$.
If  $c^1 +c^2\ge  c^3$ there is a range of   choices:

If $c^1\ge c^2$, choose $C^2\in
[c^3-c^1, c^2]$ and $C^1 =c^3 -C^2$.

If $c^1\le c^2$, choose $C^1 \in
[c^3-c^2, c^1]$ and $C^2 =c^3 -C^1$.

In the fair merging case considered above, we have   for Case C that either 
 $C^1$ is equal to the upper bound $c^1$ or $C^2$ is equal to the upper bound $c^2$. 
 For  Case B we have that $C^1=C^2 = \frac{c^3}{2}$, which is inside the above intervals. 
In case of a priority lane we have   $C^1 = c^3 \le c^1$ in Case B and $C^2 =c^3-c^1$ in Case C.
The other models might be derived from suitable kinetic coupling conditions in a similar way.
\end{remark}

\section{Numerical results}
\label{Numerical results}
In this section we test kinetic and macroscopic network solutions with the respective coupling conditions for several characteristic numerical examples. 

The two velocity model is discretized in its conservative form \eqref{eq:lindeg+relax} using a Godunov scheme, as well as the LWR model. 
In all numerical examples the intervals on the edges $[0,1]$ are discretized with $1000$ cell.
Edges $1$ and $2$ are connected to the junction at $x=1$, while the cars enter at $x=0$ into road $3$.
At the outer boundaries zero-Neumann boundary conditions are imposed.
The scaling parameter $\varepsilon $ in the kinetic equation is chosen as $\varepsilon = 0.001$.
As initial conditions the densities $\rho^i$ are chosen constant on each road.
The additional initial condition for $Z$ in the kinetic model is chosen as $Z=\frac{F(\rho^i)}{1-\rho^i+F(\rho^i)}$.  
All solutions are computed up to $T=1$.

\subsection{Fair merging}
\begin{figure}[h]
	\externaltikz{merge_case1}{
		\begin{tikzpicture}[scale=0.65]
			\begin{groupplot}[
				group style={group size=2 by 3, vertical sep = 0.75cm, horizontal sep = 1.75cm},
				width = 9cm,
				height = 4cm,
				xmin = -0.0, xmax = 1.0,
				ymin = 0.0, ymax = 1.0,
				legend style = {at={(0.5,1)},xshift=0.2cm,yshift=0.1cm,anchor=south},
				legend columns= 3,			
			]
			
			\nextgroupplot[ ylabel = $\rho^1$]
				\addplot[color = blue,thick] file{Data/merge_Lindeg_rho_1ex1_eps0001.txt};
				\addlegendentry{kinetic}
				\addplot[color = red,thick] file{Data/merge_LWR_rho_1ex1.txt};
				\addlegendentry{LWR}
			\nextgroupplot[ 
				xmin = 0.95, xmax = 1.0,
				ymin = 0.05, ymax = 0.2]	
				\addplot[color = blue,thick] file{Data/merge_Lindeg_rho_1ex1_eps0001.txt};
				\addplot[color = red,thick] file{Data/merge_LWR_rho_1ex1.txt};
				
			\nextgroupplot[  ylabel =  $\rho^2$]
				\addplot[color = blue,thick] file{Data/merge_Lindeg_rho_2ex1_eps0001.txt};
				\addplot[color = red,thick] file{Data/merge_LWR_rho_2ex1.txt};
			\nextgroupplot[
				xmin = 0.95, xmax = 1.0,
				ymin = 0.1, ymax = 0.25]	
				\addplot[color = blue,thick] file{Data/merge_Lindeg_rho_2ex1_eps0001.txt};
				\addplot[color = red,thick] file{Data/merge_LWR_rho_2ex1.txt};
			\nextgroupplot[  ylabel = $\rho^3$]
				\addplot[color = blue,thick] file{Data/merge_Lindeg_rho_3ex1_eps0001.txt};
				\addplot[color = red,thick] file{Data/merge_LWR_rho_3ex1.txt};
			\nextgroupplot[
					xmin = 0.0, xmax = 0.05,
					ymin = 0.3, ymax = 0.4]	
				\addplot[color = blue,thick] file{Data/merge_Lindeg_rho_3ex1_eps0001.txt};
				\addplot[color = red,thick] file{Data/merge_LWR_rho_3ex1.txt};
			\end{groupplot}
		\end{tikzpicture}
	}
	\caption{Fair merging with $\rho^1 = 0.1$, $\rho^2 = 0.15$, $\rho^3 = 0.2$.  }
	\label{fig:Merge_case1}
\end{figure}
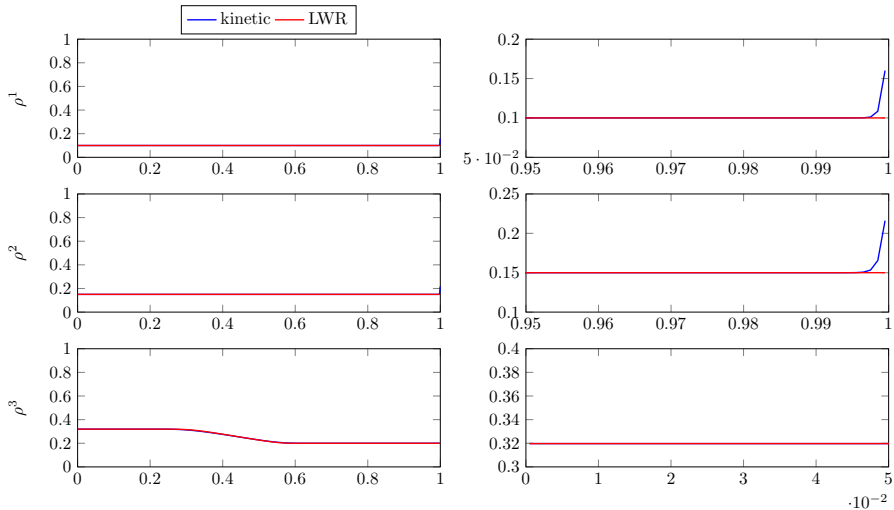

First we compare the numerical solutions of the two velocity model with the coupling conditions from section \ref{sec:fair} to the results obtained for the LWR model with the coupling conditions \eqref{eq:LWRcoupling}.
The initial densities are chosen as $\rho^1 = 0.1$, $\rho^2 = 0.15$ and  $\rho^3 = 0.2$.
The densities are small enough, such  that all  cars can pass the junction, which corresponds to Case A with $c^1+c^2 \le c^3$ and more accurately to Case 7, first subcase. The $\rho_0^i$ are given by $\rho_-(F(\rho_B^1)+F(\rho_B^2))$ with a numerical value $\rho_0^i=0.3197$.
In Figure \ref{fig:Merge_case1} the numerical solutions are shown. 
The solution of the kinetic model (blue) is almost identical to the solution of the LWR model (red).
On the right hand side there is a zoom onto the boundaries at the junction.
On edge $1$ and $2$ we can observe two boundary layers, as these correspond to stable cases.
In edge $3$ there is no layer, since the half space solution is unstable. The solution at $x=0$ fits exactly to the analytical value.

\begin{figure}[h]
	\externaltikz{merge_case2}{		
		\begin{tikzpicture}[scale=0.65]
			\begin{groupplot}[
				group style={group size=2 by 3, vertical sep = 0.75cm, horizontal sep = 1.75cm},
				width = 9cm,
				height = 4cm,
				xmin = -0.0, xmax = 1.0,
				ymin = 0.0, ymax = 1.0,
				legend style = {at={(0.5,1)},xshift=0.2cm,yshift=0.1cm,anchor=south},
				legend columns= 3,			
				]
			\nextgroupplot[ ylabel = $\rho^1$]
				\addplot[color = blue,thick] file{Data/merge_Lindeg_rho_1ex2_eps0001.txt};
				\addlegendentry{kinetic}
				\addplot[color = red,thick] file{Data/merge_LWR_rho_1ex2.txt};
				\addlegendentry{LWR}
			\nextgroupplot[ 
					xmin = 0.95, xmax = 1.0,
					ymin = 0.7, ymax = 0.9]	
				\addplot[color = blue,thick] file{Data/merge_Lindeg_rho_1ex2_eps0001.txt};
				\addplot[color = red,thick] file{Data/merge_LWR_rho_1ex2.txt};
			\nextgroupplot[ ylabel = $\rho^2$]
				\addplot[color = blue,thick] file{Data/merge_Lindeg_rho_2ex2_eps0001.txt};
				\addplot[color = red,thick] file{Data/merge_LWR_rho_2ex2.txt};
				\nextgroupplot[ 
				xmin = 0.95, xmax = 1.0,
				ymin = 0.7, ymax = 0.9]	
				\addplot[color = blue,thick] file{Data/merge_Lindeg_rho_2ex2_eps0001.txt};
				\addplot[color = red,thick] file{Data/merge_LWR_rho_2ex2.txt};
			\nextgroupplot[  ylabel = $\rho^3$]
				\addplot[color = blue,thick] file{Data/merge_Lindeg_rho_3ex2_eps0001.txt};
				\addplot[color = red,thick] file{Data/merge_LWR_rho_3ex2.txt};
			\nextgroupplot[
				xmin = 0.0, xmax = 0.05,
				ymin = 0.4, ymax = 0.9]	
				\addplot[color = blue,thick] file{Data/merge_Lindeg_rho_3ex2_eps0001.txt};
				\addplot[color = red,thick] file{Data/merge_LWR_rho_3ex2.txt};
			\end{groupplot}
		\end{tikzpicture}
	}
	\caption{Fair merging with $\rho^1 = 0.7$, $\rho^2 = 0.6$, $\rho^3 = 0.2$.  }
	\label{fig:Merge_case2}
\end{figure}
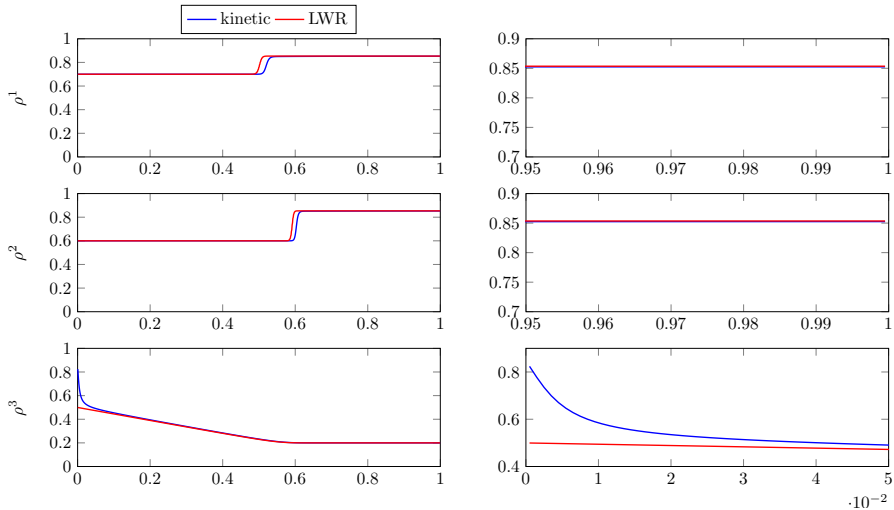

In Figure \ref{fig:Merge_case2} the numerical solutions to the initial values $\rho^1 = 0.7$, $\rho^2 = 0.6$ and  $\rho^3 = 0.2$ are shown. 
In this situation more cars are approaching the junction than can enter  road $3$. We are in the situation of Case 1 with the analytical value $\rho_0^i =\rho_+(\sigma/2)=0.83536$.
Thus the flow in the exiting road is set to its maximum, while there are jams propagating upstream in the ingoing roads.
Here we observe only in edge $3$ a layer, which interacts with the tail of the rarefaction wave.
In the ingoing roads the unstable layer solution enforce the new values at the junction.
In these roads the shock waves of the kinetic model are slightly behind those of the macroscopic one. 
This stems from an initial layer, as the layer at the junction has to form at the beginning. 
This happens in short time and is not visible at the rarefaction waves, but it remains noticeable at the shocks.
The speeds of the shocks is identical in both models, as the connected states coincide, i.e. the delay does not change over time.

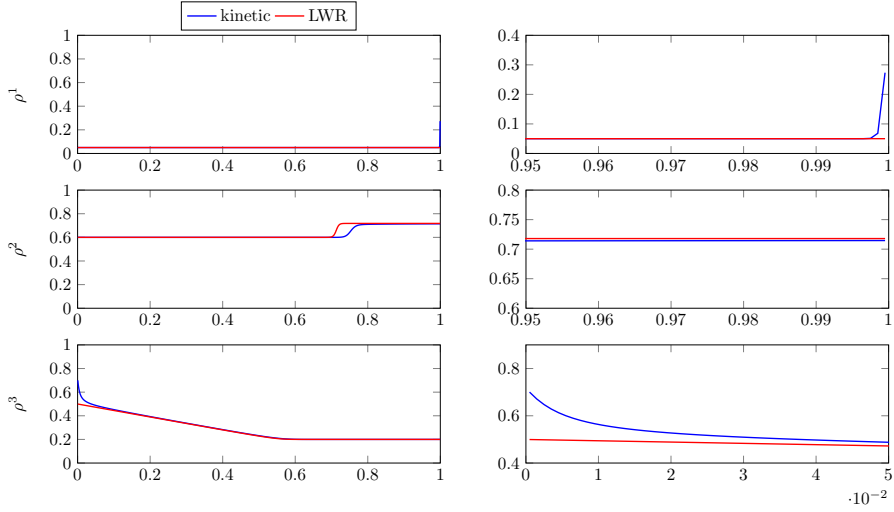
\begin{figure}[h]
	\externaltikz{merge_case3}{
		\begin{tikzpicture}[scale=0.65]
			\begin{groupplot}[
				group style={group size=2 by 3, vertical sep = 0.75cm, horizontal sep = 1.75cm},
				width = 9cm,
				height = 4cm,
				xmin = -0.0, xmax = 1.0,
				ymin = 0.0, ymax = 1.0,
				legend style = {at={(0.5,1)},xshift=0.2cm,yshift=0.1cm,anchor=south},
				legend columns= 3,			
				]		
			\nextgroupplot[ ylabel = $\rho^1$]
				\addplot[color = blue,thick] file{Data/merge_Lindeg_rho_1ex3_eps0001.txt};
				\addlegendentry{kinetic}
				\addplot[color = red,thick] file{Data/merge_LWR_rho_1ex3.txt};
				\addlegendentry{LWR}
			\nextgroupplot[ 
				xmin = 0.95, xmax = 1.0,
				ymin = 0.0, ymax = 0.4]	
				\addplot[color = blue,thick] file{Data/merge_Lindeg_rho_1ex3_eps0001.txt};
				\addplot[color = red,thick] file{Data/merge_LWR_rho_1ex3.txt};
			\nextgroupplot[ ylabel = $\rho^2$]
				\addplot[color = blue,thick] file{Data/merge_Lindeg_rho_2ex3_eps0001.txt};
				\addplot[color = red,thick] file{Data/merge_LWR_rho_2ex3.txt};
			\nextgroupplot[ 
				xmin = 0.95, xmax = 1.0,
				ymin = 0.6, ymax = 0.8]	
				\addplot[color = blue,thick] file{Data/merge_Lindeg_rho_2ex3_eps0001.txt};
				\addplot[color = red,thick] file{Data/merge_LWR_rho_2ex3.txt};
			\nextgroupplot[  ylabel = $\rho^3$]
				\addplot[color = blue,thick] file{Data/merge_Lindeg_rho_3ex3_eps0001.txt};
				\addplot[color = red,thick] file{Data/merge_LWR_rho_3ex3.txt};
			\nextgroupplot[
				xmin = 0.0, xmax = 0.05,
				ymin = 0.4, ymax = 0.9]	
				\addplot[color = blue,thick] file{Data/merge_Lindeg_rho_3ex3_eps0001.txt};
				\addplot[color = red,thick] file{Data/merge_LWR_rho_3ex3.txt};
			\end{groupplot}
		\end{tikzpicture}
	}
	\caption{Fair merging with $\rho^1 = 0.05$, $\rho^2 = 0.6$, $\rho^3 = 0.2$. }
	\label{fig:Merge_case3}
\end{figure}

In the next example, with the initial values $\rho^1 = 0.05$, $\rho^2 = 0.6$ and $\rho^3 = 0.2$, few cars enter from road $1$ but many from road $2$. We are in Case 4, first subcase. The analytical value at the junction is $\rho_0^i=\rho_+(\sigma-F(\rho_B^1))=0.7179$.
As shown in Figure \ref{fig:Merge_case3}, the flow in road $3$ is at maximum such that all cars from road $1$ and most of road $2$ can pass.
Those which do not fit in, create a jam in road $2$.
Again we see a delay of the shock, as in the previous example. 
Similarly we observe a layer in edge $3$.
But here also a layer in road $1$ is present, as the solution of the half space is now stable.

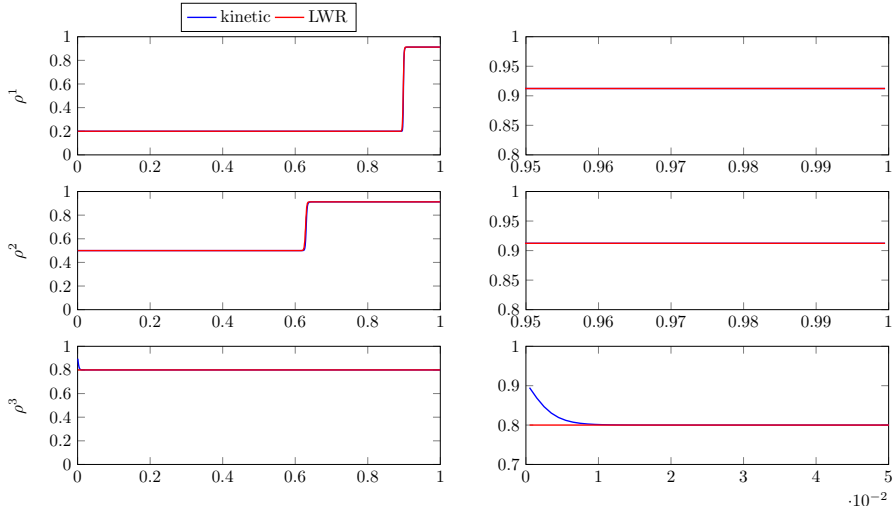
\begin{figure}[h]
	\externaltikz{merge_case31}{
		\begin{tikzpicture}[scale=0.65]
		\begin{groupplot}[
		group style={group size=2 by 3, vertical sep = 0.75cm, horizontal sep = 1.75cm},
		width = 9cm,
		height = 4cm,
		xmin = -0.0, xmax = 1.0,
		ymin = 0.0, ymax = 1.0,
		legend style = {at={(0.5,1)},xshift=0.2cm,yshift=0.1cm,anchor=south},
		legend columns= 3,			
		]		
		\nextgroupplot[ ylabel = $\rho^1$]
		\addplot[color = blue,thick] file{Data/merge_Lindeg_rho_1ex31_eps0001.txt};
		\addlegendentry{kinetic}
		\addplot[color = red,thick] file{Data/merge_LWR_rho_1ex31.txt};
		\addlegendentry{LWR}
		\nextgroupplot[ 
		xmin = 0.95, xmax = 1.0,
		ymin = 0.8, ymax = 1.0]	
		\addplot[color = blue,thick] file{Data/merge_Lindeg_rho_1ex31_eps0001.txt};
		\addplot[color = red,thick] file{Data/merge_LWR_rho_1ex31.txt};
		\nextgroupplot[ ylabel = $\rho^2$]
		\addplot[color = blue,thick] file{Data/merge_Lindeg_rho_2ex31_eps0001.txt};
		\addplot[color = red,thick] file{Data/merge_LWR_rho_2ex31.txt};
		\nextgroupplot[ 
		xmin = 0.95, xmax = 1.0,
		ymin = 0.8, ymax = 1.0]	
		\addplot[color = blue,thick] file{Data/merge_Lindeg_rho_2ex31_eps0001.txt};
		\addplot[color = red,thick] file{Data/merge_LWR_rho_2ex31.txt};
		\nextgroupplot[  ylabel = $\rho^3$]
		\addplot[color = blue,thick] file{Data/merge_Lindeg_rho_3ex31_eps0001.txt};
		\addplot[color = red,thick] file{Data/merge_LWR_rho_3ex31.txt};
		\nextgroupplot[
		xmin = 0.0, xmax = 0.05,
		ymin = 0.7, ymax = 1.0]	
		\addplot[color = blue,thick] file{Data/merge_Lindeg_rho_3ex31_eps0001.txt};
		\addplot[color = red,thick] file{Data/merge_LWR_rho_3ex31.txt};
		\end{groupplot}
		\end{tikzpicture}
	}
	\caption{Fair merging with $\rho^1 = 0.2$, $\rho^2 = 0.5$, $\rho^3 = 0.8$. }
	\label{fig:Merge_case31}
\end{figure}

If there is little space on the outgoing road, as for $\rho^1 = 0.2$, $\rho^2 = 0.5$ and $\rho^3 = 0.8$ in Figure \ref{fig:Merge_case31}, two shocks form in the ingoing roads. We are in the situation of Case 6, first subcase or Case 8, first subcase. The value at the junction is
$\rho_0^i= \rho_+(F(\rho_B^3)/2)=0.9123$.
Here the kinetic and the macroscopic solution coincide, while a layer forms on road $3$.

\subsection{Priority lane}
Now the numerical solutions of the two velocity model with the coupling conditions of section \ref{sec:prio} are compared to those obtained for the LWR model with the coupling conditions of section \ref{sec:macro_coupling_priority}.
The parameter $\delta$ is chosen $\delta=\frac{1}{2}$.
Numerical tests identical to the following ones, but with smaller values of $\delta$ show a similar behavior. 
Only the size of some layers increases, but the solution along the roads follows the same macroscopic dynamics.

In the first example with $\rho^1 = 0.6$, $\rho^2 = 0.7$ and $\rho^3 = 0.2$, shown in Figure \ref{fig:Merge_case10}, many cars arrive at the junction. 
\begin{figure}[h]
	\externaltikz{merge_case10}{
		\begin{tikzpicture}[scale=0.65]
		\begin{groupplot}[
		group style={group size=2 by 3, vertical sep = 0.75cm, horizontal sep = 1.75cm},
		width = 9cm,
		height = 4cm,
		xmin = -0.0, xmax = 1.0,
		ymin = 0.0, ymax = 1.0,
		legend style = {at={(0.5,1)},xshift=0.2cm,yshift=0.1cm,anchor=south},
		legend columns= 3,			
		]		
		\nextgroupplot[ ylabel = $\rho^1$]
		\addplot[color = blue,thick] file{Data/merge_Lindeg_rho_1ex10_eps0001.txt};
		\addlegendentry{kinetic}
		\addplot[color = red,thick] file{Data/merge_LWR_rho_1ex10.txt};
		\addlegendentry{LWR}
		\nextgroupplot[ 
		xmin = 0.95, xmax = 1.0,
		ymin = 0.2, ymax = 0.6]	
		\addplot[color = blue,thick] file{Data/merge_Lindeg_rho_1ex10_eps0001.txt};
		\addplot[color = red,thick] file{Data/merge_LWR_rho_1ex10.txt};
		\nextgroupplot[ ylabel = $\rho^2$]
		\addplot[color = blue,thick] file{Data/merge_Lindeg_rho_2ex10_eps0001.txt};
		\addplot[color = red,thick] file{Data/merge_LWR_rho_2ex10.txt};
		\nextgroupplot[ 
		xmin = 0.95, xmax = 1.0,
		ymin = 0.9, ymax = 1.1]	
		\addplot[color = blue,thick] file{Data/merge_Lindeg_rho_2ex10_eps0001.txt};
		\addplot[color = red,thick] file{Data/merge_LWR_rho_2ex10.txt};
		\nextgroupplot[  ylabel = $\rho^3$]
		\addplot[color = blue,thick] file{Data/merge_Lindeg_rho_3ex10_eps0001.txt};
		\addplot[color = red,thick] file{Data/merge_LWR_rho_3ex10.txt};
		\nextgroupplot[
		xmin = 0.0, xmax = 0.05,
		ymin = 0.45, ymax = 0.8]	
		\addplot[color = blue,thick] file{Data/merge_Lindeg_rho_3ex10_eps0001.txt};
		\addplot[color = red,thick] file{Data/merge_LWR_rho_3ex10.txt};
		\end{groupplot}
		\end{tikzpicture}
	}
	\caption{Priority merge with $\rho^1 = 0.6$, $\rho^2 = 0.7$, $\rho^3 = 0.2$. }
	\label{fig:Merge_case10}
\end{figure}
As those of road $1$ have priority, the maximal flow is established, while all cars in road $2$ have to wait. 
Layers can be observed in road $1$ and $3$. This configuration corresponds to Case B.

In the second example with $\rho^1 = 0.1$, $\rho^2 = 0.5$ and $\rho^3 = 0.2$ only few cars arrive from the first road. 
In Figure \ref{fig:Merge_case11} we can see that all the cars from road $1$ can pass, while some have to wait in the second road. 
This is a situation from Case C.
\begin{figure}[h]
	\externaltikz{merge_case11}{
		\begin{tikzpicture}[scale=0.65]
		\begin{groupplot}[
		group style={group size=2 by 3, vertical sep = 0.75cm, horizontal sep = 1.75cm},
		width = 9cm,
		height = 4cm,
		xmin = -0.0, xmax = 1.0,
		ymin = 0.0, ymax = 1.0,
		legend style = {at={(0.5,1)},xshift=0.2cm,yshift=0.1cm,anchor=south},
		legend columns= 3,			
		]		
		\nextgroupplot[ ylabel = $\rho^1$]
		\addplot[color = blue,thick] file{Data/merge_Lindeg_rho_1ex11_eps0001.txt};
		\addlegendentry{kinetic}
		\addplot[color = red,thick] file{Data/merge_LWR_rho_1ex11.txt};
		\addlegendentry{LWR}
		\nextgroupplot[ 
		xmin = 0.95, xmax = 1.0,
		ymin = 0.0, ymax = 0.2]	
		\addplot[color = blue,thick] file{Data/merge_Lindeg_rho_1ex11_eps0001.txt};
		\addplot[color = red,thick] file{Data/merge_LWR_rho_1ex11.txt};
		\nextgroupplot[ ylabel = $\rho^2$]
		\addplot[color = blue,thick] file{Data/merge_Lindeg_rho_2ex11_eps0001.txt};
		\addplot[color = red,thick] file{Data/merge_LWR_rho_2ex11.txt};
		\nextgroupplot[ 
		xmin = 0.95, xmax = 1.0,
		ymin = 0.7, ymax = 0.9]	
		\addplot[color = blue,thick] file{Data/merge_Lindeg_rho_2ex11_eps0001.txt};
		\addplot[color = red,thick] file{Data/merge_LWR_rho_2ex11.txt};
		\nextgroupplot[  ylabel = $\rho^3$]
		\addplot[color = blue,thick] file{Data/merge_Lindeg_rho_3ex11_eps0001.txt};
		\addplot[color = red,thick] file{Data/merge_LWR_rho_3ex11.txt};
		\nextgroupplot[
		xmin = 0.0, xmax = 0.05,
		ymin = 0.4, ymax = 0.9]	
		\addplot[color = blue,thick] file{Data/merge_Lindeg_rho_3ex11_eps0001.txt};
		\addplot[color = red,thick] file{Data/merge_LWR_rho_3ex11.txt};
		\end{groupplot}
		\end{tikzpicture}
	}
	\caption{Priority merge with $\rho^1 = 0.1$, $\rho^2 = 0.5$, $\rho^3 = 0.2$. }
	\label{fig:Merge_case11}
\end{figure}
As there is a free outflow in road $1$, only a layer in the outgoing road appears.

\begin{figure}[h]
	\externaltikz{merge_case12}{
		\begin{tikzpicture}[scale=0.65]
		\begin{groupplot}[
		group style={group size=2 by 3, vertical sep = 0.75cm, horizontal sep = 1.75cm},
		width = 9cm,
		height = 4cm,
		xmin = -0.0, xmax = 1.0,
		ymin = 0.0, ymax = 1.0,
		legend style = {at={(0.5,1)},xshift=0.2cm,yshift=0.1cm,anchor=south},
		legend columns= 3,			
		]		
		\nextgroupplot[ ylabel = $\rho^1$]
		\addplot[color = blue,thick] file{Data/merge_Lindeg_rho_1ex12_eps0001.txt};
		\addlegendentry{kinetic}
		\addplot[color = red,thick] file{Data/merge_LWR_rho_1ex12.txt};
		\addlegendentry{LWR}
		\nextgroupplot[ 
		xmin = 0.95, xmax = 1.0,
		ymin = 0.6, ymax = 0.8]	
		\addplot[color = blue,thick] file{Data/merge_Lindeg_rho_1ex12_eps0001.txt};
		\addplot[color = red,thick] file{Data/merge_LWR_rho_1ex12.txt};
		\nextgroupplot[ ylabel = $\rho^2$]
		\addplot[color = blue,thick] file{Data/merge_Lindeg_rho_2ex12_eps0001.txt};
		\addplot[color = red,thick] file{Data/merge_LWR_rho_2ex12.txt};
		\nextgroupplot[ 
		xmin = 0.95, xmax = 1.0,
		ymin = 0.9, ymax = 1.1]	
		\addplot[color = blue,thick] file{Data/merge_Lindeg_rho_2ex12_eps0001.txt};
		\addplot[color = red,thick] file{Data/merge_LWR_rho_2ex12.txt};
		\nextgroupplot[  ylabel = $\rho^3$]
		\addplot[color = blue,thick] file{Data/merge_Lindeg_rho_3ex12_eps0001.txt};
		\addplot[color = red,thick] file{Data/merge_LWR_rho_3ex12.txt};
		\nextgroupplot[
		xmin = 0.0, xmax = 0.05,
		ymin = 0.6, ymax = 0.9]	
		\addplot[color = blue,thick] file{Data/merge_Lindeg_rho_3ex12_eps0001.txt};
		\addplot[color = red,thick] file{Data/merge_LWR_rho_3ex12.txt};
		\end{groupplot}
		\end{tikzpicture}
	}
	\caption{Priority merge with $\rho^1 = 0.4$, $\rho^2 = 0.4$, $\rho^3 = 0.7$. }
	\label{fig:Merge_case12}
\end{figure}
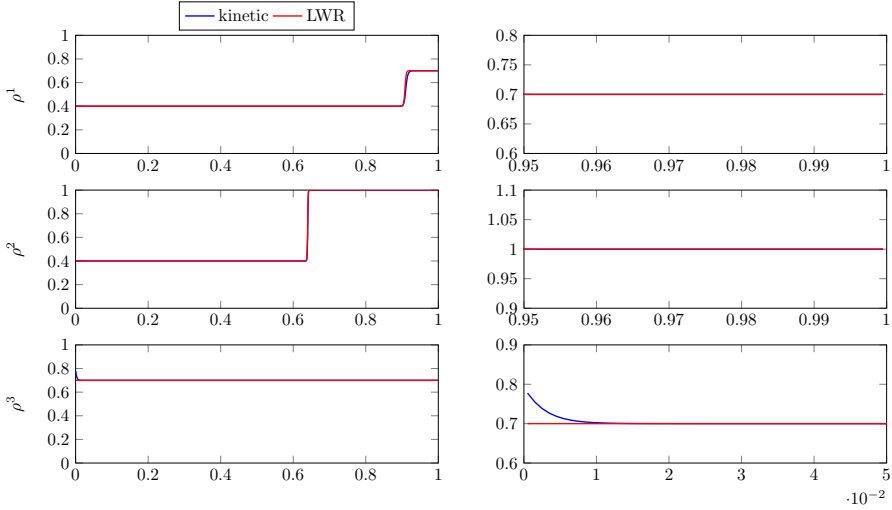

As a last example we consider a situation, where the same amount of cars is in the ingoing roads, but only little space in the outgoing one, i.e. 
$\rho^1 = 0.4$, $\rho^2 = 0.4$, $\rho^3 = 0.7$.
As expected, we can see in Figure \ref{fig:Merge_case12} that all the cars in road $2$ have to wait and thus a larger shock forms.
Not all cars in the first road can pass, but the flow is larger as in the second road. This is again a situation from Case B.

\section{Technical details}
\label{technical}

In this section we give the proof of Theorem \ref{theoremmerge}.

\subsection{Coupling the kinetic layers}		
\label{layerproof}

In the first  step of the proof of Theorem  \ref{theoremmerge}  the combination of the kinetic coupling conditions  with the layer equations has to be  considered.
The states at the junction (corresponding to $y=0$ for the layers) are denoted by
$\rho_0^i$.
Each layer can have either a stable solution (S) or an unstable solution (U).
Thus, for three edges we have eight possible combinations, which be denote by U/S-U/S-U/S.

\noindent{\bf Case1, U-U-U.}	
We  have $\rho_0^1 = \rho_+ (C^1), \rho_0^2= \rho_+(C^2) ,\rho_0^3= \rho_-(C^3)$.
The coupling conditions give
\begin{align*}
\rho_+(C^1) &= \rho_+(C^2)= \rho_-(C^3)\\
C^3&= C^1+C^2
\end{align*}	
with $0 \le  C^1,C^2,C^3< \sigma$.  

The second equality gives  $C^2= C^3=\sigma$. This is not consistent with the range of $C^2$ and $C^3$.
The case is not admissible.

\noindent{\bf Case 2, S-U-U} 
We have $\rho^1_0 \in \lceil 0,\rho_+ (C^1)\rfloor$ and $ \rho^2_0= \rho_+(C^2), \rho^3_0= \rho_-(C^3)$.
Inserting into the coupling conditions gives
\begin{eqnarray*}
	\rho_0^1 &= &\rho_+(C^2)= \rho_-(C^3)\\
	C^3 & =&C^1+C^2
\end{eqnarray*}	
with $0 \le  C^1 \le  \sigma$ and $0 \le  C^2,C^3< \sigma$.
Again the second  equation gives $C^2=C^3= \sigma$ which is not in the range of  $C^2,C^3$. 
The case is not admissible.

\noindent{\bf Case 3, U-S-U} 
We have $\rho_0^1 = \rho_+ (C^1), \rho_0^2 \in [0,\rho_+(C^2)), \rho_0^3= \rho_-(C^3)$.
The case is symmetric to the above and not admissible.

\noindent{\bf Case 4, U-U-S} 
We have $\rho_0^1 = \rho_+ (C^1), \rho_0^2 = \rho_+(C^2)$, $ \rho_0^3= \lceil\rho_-(C^3),1\rfloor$.
We have
\begin{align*}
\rho_+(C^1) &= \rho_+(C^2)= \rho_0^3\\
C^3&=C^1+C^2
\end{align*}
with $0 \le C^1,C^2 < \sigma$	 and $0 \le C^3 \le \sigma$.
This gives $C^1=C^2= \frac{C^3}{2}$ and $$	\rho_0^1 = 	\rho_0^2 =	\rho_0^3 =\rho_+(\frac{C^3}{2}).$$

\noindent{\bf Case 5, U-S-S} 
We have $\rho_0^1 = \rho_+ (C^1), \rho_0^2 \in \lceil0,\rho_+(C^2)\rfloor$ 
 and $ \rho_0^3\in  \lceil\rho_-(C^3),1\rfloor$ with 
$0 \le C^1 < \sigma$	 and $0 \le C^2, C^3 \le \sigma$.
We have 
\begin{align*}
\rho_+(C^1) &= \rho_0^2= \rho_0^3\\
C^3&=C^1+C^2.
\end{align*}	
This gives $\rho_0^2 = \rho_0^3 = \rho_+(C^1) =\rho_+(C^3-C^2)$ with
the requirement $0 \le C^3-C^2 \le \sigma$ or
$ C^3 \ge C^2 $ and $\rho_0^2 = \rho_0^3 =\rho_+(C^3-C^2) \in [\rho_-(C^3),\rho_+(C^2)]$. It leads to $\rho_+(C^3-C^2) \le \rho_+(C^2)$ or $C^3-C^2 \ge  C_2$ or $C^3 \ge  2 C^2$.
Altogether, we have for $ 2 C^2 \le C^3$  and  $C^1 = C^3-C^2$
\begin{eqnarray*}
	\rho_0^1&=\rho_0^2 =	\rho_0^3 = \rho_+(C^3-C^2)\ .
\end{eqnarray*}

\noindent{\bf Case 6, S-U-S} 
We have $\rho^1_0 \in  \lceil0,\rho_+ (C^1)\rfloor$  and $\rho_0^2 =\rho_+(C^2), \rho_0^3\in  \lceil\rho_-(C^3),1\rfloor$ with
$0 \le C^1 \le \sigma$	 and $0 \le C^2, C^3 < \sigma$.
The case is symmetric to case 5.
For $ 2 C^1 \le C^3  $ and $C^2 = C^3-C^1$ we have
\begin{eqnarray*}
	\rho_0^1=\rho_0^2 &=	\rho_0^3 = \rho_+(C^3-C^1)\ .
\end{eqnarray*}

\noindent{\bf Case 7, S-S-U}
We have $\rho_0^1 \in \lceil 0,\rho_+ (C^1)\rfloor $  and $\rho_0^2 \in \lceil0,\rho_+(C^2)\rfloor$ 
and $ \rho_0^3 =   \rho_-(C^3)$ with $0 \le C^1,C^2 \le \sigma$	 and $0 \le  C^3 < \sigma$.
The coupling conditions give 
\begin{align*}
\rho_0^1 &= \rho_0^2= \rho_-(C^3)\\
C^3&=C^1+C^2\ .	
\end{align*}	
This gives $\rho_0^1 = \rho_0^2 = \rho_-(C^3) $
with  the condition $0 \le C^1+C^2 < \sigma$.  Thus, for $0 \le C^1+C^2 < \sigma$ we have 
$$\rho_0^1 = \rho_0^2 = \rho_0^3 = \rho_-(C^1+C^2)\ .$$

\noindent{\bf Case 8, S-S-S} 
We have $\rho_0^1 \in \lceil 0,\rho_+ (C^1)\rfloor, \rho_0^2 \in \lceil0,\rho_+(C^2)\rfloor, \rho_0^3 \in\lceil   \rho_-(C^3),1\rfloor$
with $0 \le C^1,C^2,C^3 \le \sigma$.
The conditions are  
\begin{align*}
\rho_0^1 &= \rho_0^2= \rho_0^3\\
C^3&=C^1+C^2.	
\end{align*}	
The values of $\rho^1_0 = \rho^2_0=\rho^3$ are not uniquely determined, but  they restricted to  the interval $  [\rho_-(C^1+C^2),\min(\rho_+(C^1),\rho_+(C^2))]$.	

These considerations yield all possible combinations of layer problems at the node. They have now to be matched
to the half-Riemann problems at the respective lanes.

\subsection{Matching of  Riemann problem and layer equations}
\label{proof}

Assuming  the  initial states $\rho_B^i, i=1, 2,3$ to be given, we have to determine the fluxes $C^i$ and new states $\rho_K^i$ at the node. As mentioned, on the one hand $\rho_K^i$ are the asymptotic states of the respective layer problems and they and  the corresponding fluxes $C^i$ have to fulfill the conditions on the single kinetic layers, see section \ref{summary}, and on the coupled layers, see section \ref{layerproof}.
On the other hand they are the left (road 1 and 2) or right hand (road 3) states of the half Riemann problems with $\rho_B^i$
as corresponding right and left states respectively fullfillling the conditions in \ref{Riemann}.

As before, the states at the junction (corresponding to $y=0$ for the layers) are denoted by
$\rho_0^i$.
We  consider eight different configurations for the states $\rho_B^i$ corresponding to the possible combinations of different half Riemann problems.
For  each of them all possible combinations with stable or unstable layer solutions have to be discussed.
Not admissible combinations are not listed.

\noindent{\bf Case 1, RP1-1-1} $\rho_B^1 \ge \rho^\star , \rho_B^2 \ge \rho^\star , \rho_B^3 \le  \rho^\star $.
From Section \ref{summary} we obtain
\begin{align*}
\rho_K^1 &\in [\rho^\star,1] :
&  (U) &\text{ or } ((S) \text{ with } C^1=\sigma) \\
\rho_K^2 &\in [\rho^\star,1]: 
& (U) &\text{ or } ((S) \text{ with } C^2=\sigma)\\
\rho_K^3 &\in [0,\rho^\star]:
& (U) &\text{ or } ((S) \text{ with } C^3=\sigma)
\end{align*}	
Then, the discussion in Section \ref{layerproof} leads to 5 different cases:
\begin{enumerate}
	\item[{\bf UUS}] with $C^3 = \sigma$  and  $C^1=C^2=\frac{\sigma}{2}$ and $\rho_0^3 =\rho_+(\frac{\sigma}{2})$.
	\item[{\bf USS}] with $C^2 = C^3 = \sigma$ 
	which contradicts $C^3 \ge 2 C^2$.
	\item[{\bf SUS}] with $C^1 = C^3 = \sigma$ which contradicts  $C^3 \ge 2 C^1$.
	\item[{\bf SSU}] with $C^1 = C^2 = \sigma$ and a contradiction to $C^1+C^2 \le \sigma$. 
	\item[{\bf SSS}] with $C^1 = C^2 = C^3 = \sigma$, which gives a contradiction to the balance of fluxes.		
\end{enumerate}
 This and again Section \ref{summary} give
\begin{align*}
\rho_K^1 &= \rho_+(\frac{\sigma}{2}) &\rho_K^2 &=\rho_+(\frac{\sigma}{2} ) &\rho_K^3 &=  \rho_+(\sigma)=\rho^\star\\
\rho_0^i &= \rho_+(\frac{\sigma}{2})\ . &&&& 
\end{align*}

\noindent{\bf Case 2, RP1-1-2} $\rho_B^1 \ge \rho^\star , \rho_B^2 \ge \rho^\star , \rho_B^3 \ge  \rho^\star $.
\begin{align*}
\rho_K^1 &\in [\rho^\star,1] :
&  (U) &\text{ or } ((S) \text{ with } C^1=\sigma) \\
\rho_K^2 &\in [\rho^\star,1]: 
& (U) &\text{ or } ((S) \text{ with } C^2=\sigma)\\		
\rho_K^3 &\in [0,\tau(\rho_B^3)] \cup \{\rho_B^3\}:
& ((U) \text{ with } C^3 \le F(\rho_B^3)) &\text{ or } ((S)
\text{ with } C^3=F(\rho_B^3))
\end{align*}	
\begin{enumerate}
	\item[{\bf UUS}] with $C^3 = F(\rho_B^3)$  and  $C^1=C^2= \frac{1}{2} F(\rho_B^3)$ and $\rho_0^3 =\rho_+(\frac{1}{2}F(\rho_B^3)).$
	\item[{\bf USS}] with $C^2 =  \sigma$ 
	which contradicts $C^3 \ge 2 C^2$.
	\item[{\bf SUS}] with $C^1 =  \sigma$ which contradicts  $C^3 \ge 2 C^1$.
	\item[{\bf SSU}] with $C^1 = C^2 = \sigma$ and a contradiction to $C^1+C^2 \le \sigma$. 
	\item[{\bf SSS}] with $C^1 = C^2 = C^3 = \sigma$, which gives a contradiction to the balance of fluxes.		
\end{enumerate}
This gives
\begin{align*}
\rho_K^1 &= \rho_+(\frac{1}{2}F(\rho_B^3)) &\rho_K^2 &=\rho_+(\frac{1}{2}F(\rho_B^3)) &\rho_K^3 &=  \rho_+(F(\rho_B^3))=\rho_B^3\\
\rho_0^i &=  \rho_+(\frac{1}{2}F(\rho_B^3))\ .&&&& 
\end{align*}

\noindent{\bf Case 3, RP1-2-1} $\rho_B^1 \ge \rho^\star , \rho_B^2 \le \rho^\star , \rho_B^3 \le  \rho^\star $. 
\begin{align*}
\rho_K^1 &\in [\rho^\star,1] :
&  (U) &\text{ or } ((S) \text{ with } C^1=\sigma) \\
\rho_K^2 &\in  [\tau(\rho_B^2),1] \cup \{\rho_B^2\}: 
& ((U) \text{ with } C^2 \le F(\rho_B^2)) &\text{ or } ((S)
\text{ with } C^2=F(\rho_B^2))	\\	
\rho_K^3 &\in [\rho^\star,1]:
& (U) &\text{ or } ((S) \text{ with } C^3=\sigma)	
\end{align*}	
\begin{enumerate}
	\item[{\bf UUS}] with $C^3 = \sigma$ which gives  $C^1 = C^2= \frac{\sigma}{2}$. Moreover,  $C^2\le F(\rho_B^2)$. This is possible, if $\frac{\sigma}{2} \le F(\rho_B^2)$.
	Then, $\rho_0^3 = \rho_+(\frac{\sigma}{2})$
	\item[{\bf USS}] with $C^3 =  \sigma, C^2=F(\rho_B^2)$. 
	$C^3 \ge 2 C^2$ gives the requirement  $\frac{\sigma}{2} \ge F(\rho_B^2)$. Moreover, we have $C^1 = \sigma-F(\rho_B^2)  $ and 
	$\rho_0^2 = \rho_0^3= \rho_+(C^1)$.
	\item[{\bf SUS}] with $C^1 =  \sigma$ which contradicts  $C^3 \ge 2 C^1$.
	\item[{\bf SSU}] with $C^1 =  \sigma$ and $C^2 = F(\rho_B^2)$. This is only possible for $\rho_B^2 = 0$
	and $C^2 =0$. Then $C^3 =\sigma$ and  $\rho_0^i = \rho_-(\sigma) = \rho^\star$.
	\item[{\bf SSS}] with $C^1 =  C^3 = \sigma$ and $C^2 = F(\rho_B^2)$. This gives again $\rho_B^2=0$ and
	$C^2=0$. Then  $\rho_0^i \in [\rho_-(C^1+C^2),\min(\rho_+(C^1),\rho_+(C^2))]$ gives 
	$\rho_0^i \in [\rho_-(\sigma),\rho_+(\sigma)]$. This  leaves only $\rho_0^i =\rho^\star$.
\end{enumerate}
This gives  for $\frac{\sigma}{2} \le F(\rho_B^2) $ that   $C^1 = C^2= \frac{C^3}{2}=\frac{\sigma}{2}$ and 
\begin{align*}
\rho_K^1 &= \rho_+(\frac{\sigma}{2}) &\rho_K^2 &=\rho_+(\frac{\sigma}{2}) &\rho_K^3 &=  \rho_+(\sigma)=\rho^\star\\
\rho_0^i &=  \rho_+(\frac{\sigma}{2})\ .&&&& 
\end{align*}

For 	$\frac{\sigma}{2} \ge F(\rho_B^2) $		one has	$C^1 = \sigma-F(\rho_B^2)  $ , $C^3 =  \sigma, C^2=F(\rho_B^2)$ and 
\begin{align*}
\rho_K^1 &= \rho_+(\sigma-F(\rho_B^2)) &\rho_K^2 &=\rho_-(F(\rho_B^2)) =\rho_B^2&\rho_K^3 &=  \rho_+(\sigma)=\rho^\star\\
\rho_0^i &=  \rho_+( \sigma-F(\rho_B^2) )\ .&&&& 
\end{align*}

\noindent{\bf Case 4, RP2-1-1} $\rho_B^1 \le \rho^\star , \rho_B^2 \ge \rho^\star , \rho_B^3 \le  \rho^\star $. 
This case is symmetric to Case 3.

We have  for $\frac{\sigma}{2} \ge F(\rho_B^1)$ that $C^1 = F(\rho_B^1) , C^3 =\sigma$ , $C^2 = \sigma - F(\rho_B^1)$ and 
\begin{align*}
\rho_K^1 &= \rho_B^1 &\rho_K^2 &=\rho_+(\sigma-F(\rho_B^1)) &\rho_K^3 &=  \rho_+(\sigma)=\rho^\star\\
\rho_0^i &=  \rho_+(\sigma-F(\rho_B^1))\ .&&&&
\end{align*}

For 	$\frac{\sigma}{2} \le F(\rho_B^1) $		one has	$C^1 = C^2= \frac{C^3}{2} = \frac{\sigma}{2}$ and 
\begin{align*}
\rho_K^1 &= \rho_+(\frac{\sigma}{2}) &\rho_K^2 &=\rho_+(\frac{\sigma}{2}) &\rho_K^3 &=  \rho_+(\sigma)=\rho^\star\\
\rho_0^i &=  \rho_+(\frac{\sigma}{2})\ .&&&& 
\end{align*}

\noindent{\bf Case 5, RP1-2-2} $\rho_B^1 \ge \rho^\star , \rho_B^2 \le \rho^\star , \rho_B^3 \ge  \rho^\star .$
\begin{align*}
\rho_K^1 &\in  [\rho^\star,1] : 
& (U)  &\text{ or } ((S)
\text{ with } C^1=\sigma)	\\
\rho_K^2 &\in [\tau(\rho_B^2),1] \cup \{\rho_B^2\}:  
&  ((U)\text{ with } C^2 \le F(\rho_B^2)) &\text{ or } ((S) \text{ with } C^2= F(\rho_B^2))\\
\rho_K^3 &\in [0,1-\rho_B^3] \cup \{\rho_B^3\}: 
&((U)\text{ with } C^3 \le F(\rho_B^3)) &\text{ or } ((S) \text{ with } C^3= F(\rho_B^3))
\end{align*}	
\begin{enumerate}
	\item[{\bf UUS}] with $C^2 \le F(\rho_B^2) $ and $C^3= F(\rho_B^3)$. If  $ F(\rho_B^3)\le 2  F(\rho_B^2)$ then
	$C^1= C^2 = \frac{F(\rho_B^3)}{2}$ and 
	$\rho_0^3 = \rho_+(\frac{C^3}{2})$.
	\item[{\bf USS}] with $C^2 = F(\rho_B^2) , C^3 = F(\rho_B^3) $. With $C^3\ge 2 C^2$ or $F(\rho_B^3) \ge F(\rho_B^2)$ we have $C^1 =C^3-C^2$. 
	\item[{\bf SUS}] with $C^1 = \sigma , C^2 \le F(\rho_B^2), C^3 = F(\rho_B^3)$, which gives  a contradiction to  $C^3 \ge 2 C^1$.
	\item[{\bf SSU}] with $C^1  = \sigma$ and $C^2 = F(\rho_B^2),C^3 \le  F(\rho_B^3) $. This is only possible for $\rho_B^2 = 0$. Then $C^3 =\sigma, \rho_B^3 = \rho^\star$ and $\rho_0^1=\rho_-(\sigma) = \rho^\star$.
	\item[{\bf SSS}] with $C^1 = \sigma$ and $C^2 =F(\rho_B^2) , C^3 =F(\rho_B^3) $. This is only possible, if $C^2 =0$ and $\rho_B^2 =0$. This yields $C^3=\sigma$ and $\rho_B^3= \rho^\star$.
	Then  $\rho_0^i \in [\rho_-(C^1+C^2),\min(\rho_+(C^1),\rho_+(C^2))]$ gives 
	$\rho_0^i \in [\rho_-(\sigma),\rho_+(\sigma)]$,	which leaves only $\rho_0^i =\rho^\star$.
\end{enumerate}

This gives  for $F(\rho_B^3)  \le  2 F(\rho_B^2)$
\begin{align*}
\rho_K^1 &= \rho_+(\frac{F(\rho_B^3)}{2}) &\rho_K^2 &=\rho_+(\frac{F(\rho_B^3)}{2} ) &\rho_K^3 &=  \rho_+(F(\rho_B^3))=\rho_B^3\\
\rho_0^i &=  \rho_+(\frac{F(\rho_B^3)}{2})\ .&&&& 
\end{align*}

For $F(\rho_B^3)  \ge  2 F(\rho_B^2) $	one has	
\begin{align*}
\rho_K^1 &= \rho_+(F(\rho_B^3)-F(\rho_B^2)) &\rho_K^2 &=\rho_-(F(\rho_B^2)) =\rho_B^2&\rho_K^3 &=  \rho_+(F(\rho_B^3))=\rho_B^3\\
\rho_0^i &=  \rho_+(F(\rho_B^3)- F(\rho_B^2))\ .&&&& 
\end{align*}

\noindent{\bf Case 6, RP2-1-2} $\rho_B^1 \le \rho^\star , \rho_B^2 \ge \rho^\star , \rho_B^3 \ge  \rho^\star $.
This case is symmetric to case 5.

We have for $F(\rho_B^3)  \le  2 F(\rho_B^1)$
\begin{align*}
\rho_K^1 &= \rho_+(\frac{F(\rho_B^3)}{2}) &\rho_K^2 &=\rho_+(\frac{F(\rho_B^3)}{2} ) &\rho_K^3 &=  \rho_+(F(\rho_B^3))=\rho_B^3\\
\rho_0^i &=  \rho_+(\frac{F(\rho_B^3)}{2})\ .&&&& 
\end{align*}

For $F(\rho_B^3 ) \ge  2 F(\rho_B^1) $	one has	
\begin{align*}
\rho_K^1 &= \rho_-(F(\rho_B^1)) =\rho_B^1&\rho_K^2 &=\rho_+(F(\rho_B^3)- F(\rho_B^1)) &\rho_K^3 &=  \rho_+(F(\rho_B^3))= \rho_B^3\\
\rho_0^i &=  \rho_+(F(\rho_B^3)- F(\rho_B^1))\ .&&&& 
\end{align*}

\noindent{\bf Case 7, RP2-2-1}  $\rho_B^1 \le \rho^\star , \rho_B^2 \le \rho^\star , \rho_B^3 \le  \rho^\star $. 
\begin{align*}
\rho_K^1 &\in [\tau(\rho_B^1),1] \cup \{\rho_B^1\}:  
&  ((U)\text{ with } C^1 \le F(\rho_B^1)) &\text{ or } ((S) \text{ with } C^1= F(\rho_B^1))\\
\rho_K^2 &\in  [\tau(\rho_B^2),1] \cup \{\rho_B^2\}: 
&((U)\text{ with } C^2 \le F(\rho_B^2)) &\text{ or } ((S) \text{ with } C^2= F(\rho_B^2))\\
\rho_K^3 &\in [0,\rho^\star] : 
& (U)  &\text{ or } ((S)
\text{ with } C^3=\sigma)	 
\end{align*}	
\begin{enumerate}
	\item[{\bf UUS}] with $C^1 \le F(\rho_B^1) $ and $C^2 \le  F(\rho_B^2) $.  $C^3=\sigma$ yields 
	$C^1= C^2 = \frac{\sigma}{2}$, if $F(\rho_B^1) \ge \frac{\sigma}{2}$ and $F(\rho_B^2) \ge \frac{\sigma}{2}$ .
	Then $\rho_0^3 = \rho_+(\frac{C^3}{2})$.
	\item[{\bf USS}] with $C^1 \le  F(\rho_B^1) , C^3 =\sigma, C^2 =  F(\rho_B^2) $.  $C^3\ge 2 C^2$ is equivalent to
	$\frac{\sigma}{2} \ge F(\rho_B^2)$. Moreover, $C^1 = \sigma - F(\rho_B^2) $ requires $F(\rho_B^1) +F(\rho_B^2) \ge \sigma$.
	\item[{\bf SUS}]  with $C^1 =  F(\rho_B^1), C^2 \le  F(\rho_B^2), C^3 =\sigma$.  $C^3 \ge 2 C^1$ gives $ \frac{\sigma}{2} \ge F(\rho_B^1)$,  $F(\rho_B^2) \ge \frac{\sigma}{2}$ and $C^2 =\sigma - F(\rho_B^1) \ge \frac{\sigma}{2}$.
	Moreover, $\rho_0^1= \rho_+(C^3-C^1)$.
	\item[{\bf SSU}] with $C^1  = F(\rho_B^1)$ and $C^2 =  F(\rho_B^2) $.  This  gives $F(\rho_B^1)+ F(\rho_B^2)\le \sigma$ and  $\rho_0^i =
	\rho_-(C^3) $.
	\item[{\bf SSS}] with $C^1 = F(\rho_B^1) $ and $C^2 =F(\rho_B^2) , C^3 =\sigma$. This is only possible, if $ F(\rho_B^1)+ F(\rho_B^2)= \sigma$. In this case, since $ \rho_0^i\in  [\rho_-(C^1+C^2),\min(\rho_+(C^1),\rho_+(C^2))]$ we obtain
	$ \rho_0^i\in  [\rho_-(\sigma),\min(\rho_+(F(\rho_B^1)),\rho_+(F(\rho_B^2)))]$. This gives the restriction
	$ \rho_0^i\in $$ [\rho^\star,\min(\tau(\rho_B^1),\tau(\rho_B^2))]$ according to the range of $\rho_B^1, \rho_B^2 $.
\end{enumerate}
We obtain for $F(\rho_B^1)+F(\rho_B^2) \le \sigma (SSU)$ 
\begin{align*}
\rho_K^1 &= \rho_-(F(\rho_B^1) )=\rho_B^1&\rho_K^2 &=\rho_-(F(\rho_B^2)) =\rho_B^2&\rho_K^3 &=  \rho_-(F(\rho_B^1)+ F(\rho_B^2))\\
\rho_0^i &=  \rho_-(F(\rho_B^1)+ F(\rho_B^2)).&&&& \ .
\end{align*}
For $F(\rho_B^1)+F(\rho_B^2) \ge \sigma ,F(\rho_B^1) \ge \frac{\sigma}{2}, F(\rho_B^2) \ge \frac{\sigma}{2} (UUS)$	one has	
\begin{align*}
\rho_K^1 &= \rho_+(\frac{\sigma}{2}) &\rho_K^2 &=\rho_+(\frac{\sigma}{2})&\rho_K^3 &=  \rho^\star\\
\rho_0^i &=  \rho_+(\frac{\sigma}{2}).&&&& \ .
\end{align*}
For $F(\rho_B^1)+F(\rho_B^2) \ge \sigma ,F(\rho_B^1) \le \frac{\sigma}{2}, F(\rho_B^2) \ge \frac{\sigma}{2} (SUS)$	one has	
\begin{align*}
\rho_K^1 &= \rho_-(F(\rho_B^1)) =\rho_B^1&\rho_K^2 &=\rho_+(\sigma-F(\rho_B^1))&\rho_K^3 &=  \rho_+(\sigma)=\rho^\star\\
\rho_0^i &=  \rho_+(\sigma-F(\rho_B^1)).&&&& \ .
\end{align*}
For $F(\rho_B^1)+F(\rho_B^2) \ge \sigma ,F(\rho_B^1) \ge \frac{\sigma}{2}, F(\rho_B^2) \le \frac{\sigma}{2}(USS)$	one has	
\begin{align*}
\rho_K^1 &= \rho_+(\sigma-F(\rho_B^2)) &\rho_K^2 &=\rho_-( F(\rho_B^2)) = \rho_B^2&\rho_K^3 &=  \rho^\star\\
\rho_0^i &=  \rho_+(\sigma- F(\rho_B^2)).&&&& \ .
\end{align*}		

\begin{remark}
We note that at the interfaces between the different conditions we obtain values  $\rho_0^i  \in   [\rho^\star,\min(\rho_+(F(\rho_B^1)),\rho_+(F(\rho_B^2))]$. This is exactly   the interval for the $\rho^i$-values in case (SSS).
\end{remark}

\noindent{\bf Case 8, RP2-2-2} $\rho_B^1 \le \rho^\star , \rho_B^2 \le \rho^\star , \rho_B^3 \ge  \rho^\star $.
\begin{align*}
\rho_K^1 &\in [\tau(\rho_B^1),1] \cup \{\rho_B^1\}:  
&  ((U)\text{ with } C^1 \le F(\rho_B^1)) &\text{ or } ((S) \text{ with } C^1= F(\rho_B^1))\\
\rho_K^2 &\in [\tau(\rho_B^2),1] \cup \{\rho_B^2\}:  
&  ((U)\text{ with } C^2 \le F(\rho_B^2)) &\text{ or } ((S) \text{ with } C^2= F(\rho_B^2))\\
\rho_K^3 &\in [0,\tau(\rho_B^3)] \cup \{\rho_B^3\}: 
&((U)\text{ with } C^3 \le F(\rho_B^3)) &\text{ or } ((S) \text{ with } C^3= F(\rho_B^3))
\end{align*}	
\begin{enumerate}
	\item[{\bf UUS}] with $C^1 \le F(\rho_B^1) $,  $C^2 \le F(\rho_B^2) $ and $C^3 = F(\rho_B^3) $.  If  $ F(\rho_B^3)\le 2  F(\rho_B^1)$ and $ F(\rho_B^3)\le 2  F(\rho_B^2)$  then
	$C^1= C^2 = \frac{F(\rho_B^3)}{2}$ and 
	$\rho_0^3 = \rho_+(\frac{C^3}{2})$.
	\item[{\bf USS}] with $C^1 \le  F(\rho_B^1) , C^2 =F(\rho_B^2), C^3 = F(\rho_B^3) $. With $C^3\ge 2 C^2$ we have 
	$ F(\rho_B^3)\ge 2  F(\rho_B^2)$ and $ F(\rho_B^3)-  F(\rho_B^2) \le F(\rho_B^1)$ 
	or $ F(\rho_B^1)+ F(\rho_B^2) \ge F(\rho_B^3)$. 
	\item[{\bf SUS}] with $C^1 =  F(\rho_B^1), C^2\le  F(\rho_B^2),C^3 = F(\rho_B^3)$.  $C^3 \ge 2 C^1$ gives $ F(\rho_B^3)\ge 2  F(\rho_B^1)$ and $F(\rho_B^3)- F(\rho_B^1) \le F(\rho_B^2)$ or  $F(\rho_B^1)+ F(\rho_B^2) \ge F(\rho_B^3)$. Moreover $\rho_0^1= \rho_+(C^3-C^1)$.
	\item[{\bf SSU}] with $C^1  = F(\rho_B^1)$ and $C^2 =  F(\rho_B^2) ,C^3 \le  F(\rho_B^3) $. 
	This is only possible for $F(\rho_B^1)+F(\rho_B^2) \le F(\rho_B^3)$. Moreover, $\rho_0^i =\rho_-(C^3) $.
	\item[{\bf SSS}] with $C^1 = F(\rho_B^1) $ and $C^2 =F(\rho_B^2), C^3 =F(\rho_B^3) $. This is only possible, if $F(\rho_B^1)+F(\rho_B^2) = F(\rho_B^3)$.
	
We obtain
	$ \rho_0^i\in  [\rho_-(F(\rho_B^3)),\min(\rho_+(F(\rho_B^1)),\rho_+(F(\rho_B^2)))]$. This gives 
	according to the range of $\rho_B^2, \rho_B^3 $, that $ \rho_0^i\in  [\tau(\rho_B^3),\min(\tau(\rho_B^1),\tau(\rho_B^2))]$.
\end{enumerate}
This gives  for $ F(\rho_B^3)\le 2  F(\rho_B^1)$ and $ F(\rho_B^3)\le 2  F(\rho_B^2) (UUS)$
\begin{align*}
\rho_K^1 &= \rho_+(\frac{F(\rho_B^3)}{2}) &\rho_K^2 &=\rho_+(\frac{F(\rho_B^3)}{2} ) &\rho_K^3 &=  \rho_+(F(\rho_B^3))=\rho_B^3\\
\rho_0^i &=  \rho_+(\frac{F(\rho_B^3)}{2})\ .&&&& 
\end{align*}

For $ F(\rho_B^3)\ge 2  F(\rho_B^2)$ and  $ F(\rho_B^1)+ F(\rho_B^2) \ge F(\rho_B^3)(USS)$	one has	
\begin{align*}
\rho_K^1 &= \rho_+(F(\rho_B^3)- F(\rho_B^2)) &\rho_K^2 &=\rho_-(F(\rho_B^2))&\rho_K^3 &=  \rho_+(F(\rho_B^3))\\
\rho_0^i &=  \rho_+(F(\rho_B^3)- F(\rho_B^2))\ .&&&& 
\end{align*}

For $F(\rho_B^3)\ge 2  F(\rho_B^1)$ and  $F(\rho_B^1)+ F(\rho_B^2) \ge F(\rho_B^3) (SUS)$	one has	
\begin{align*}
\rho_K^1 &= \rho_-(F(\rho_B^1))= \rho_B^1&\rho_K^2 &=\rho_+(F(\rho_B^3)- F(\rho_B^1)) &\rho_K^3 &=  \rho_+(F(\rho_B^3))=\rho_B^3\\
\rho_0^i &=  \rho_+(F(\rho_B^3)- F(\rho_B^1))\ .&&&& 
\end{align*}						

For  $F(\rho_B^1)+F(\rho_B^2) \le F(\rho_B^3)(SSU)$	one has	
\begin{align*}
\rho_K^1 &= \rho_-(F(\rho_B^1)) =\rho_B^1&\rho_K^2 &=\rho_-(F(\rho_B^2))=\rho_B^2 &\rho_K^3 &=  \rho_-(F(\rho_B^1)+F(\rho_B^2))\\
\rho_0^i &=  \rho_-(F(\rho_B^1)+F(\rho_B^2))\ .&&&& 
\end{align*}	

\begin{remark}
Note that the sub-cases in Case 8  partition uniquely the range of admissible states since for $0 \le x,y,z \le 1$ either ($x+y\le z$) or ($x+y\ge z$ and $z\ge2y$) or ($x+y\ge z$ and $z\ge 2x$) or ($z\le 2x$ and $z\le 2y$).		

Moreover, note that at the interfaces between the different conditions we obtain that $\rho_0^i \in [\rho_-(F(\rho_B^3)),
\min(\rho_+(F(\rho_B^1)),\rho_+(F(\rho_B^2)))]$. This is exactly   the interval for the $\rho^i$-values in case (SSS). 
\end{remark}

\section{Conclusions}
	We have introduced  coupling conditions for a  kinetic  two velocity traffic model, which is used as a relaxation model
	for scalar traffic flow equations.
	From these coupling conditions we have derived, via  asymptotic analysis of the spatial layers at the nodes and a detailed investigation of the associated Riemann problems,  coupling conditions for classical scalar macroscopic traffic models.	
	The  derivation shows that a classical merge condition  for a nonlinear scalar conservation law
	can be interpreted on the kinetic level as a combination of the balance of fluxes and the equality of densites on all roads.
The equality of densities is, however, not fulfilled on the  level of the scalar conservation law, where 
	only the balance of fluxes still holds.
	This research will be continued in a follow-up paper \cite{BK20}, where the case of diverging lanes with and whithout drivers preferences 
	is treated.


\begin{thebibliography}{10}
	\bibliographystyle{siam}
		
	\bibitem{AR}
	A.~Aw and M.~Rascle,  {\em Resurrection of second order models of traffic flow?}, SIAM J. Appl. Math., 60, 916--938, 2000.

	\bibitem{AKMR}
	 A.~Aw, A. Klar, T. Materne,  M.~Rascle, {\em Derivation of continuum flow traffic models from microscopic Follow the leader models}, SIAM J. Appl. Math. 63 (1), 259-278, 2002 

	\bibitem{AM04} 
	D. Aregba-Driollet,V. Milisic,
	{\em Kinetic approximation of a boundary value problem for conservation laws}, Numer. Math.  97, 595–633, 2004

	\bibitem{BHK06a}
	M. Banda, M. Herty, A. Klar, {\em Gas flow in pipeline networks}, NHM 1(1), 41-56, 2006	

	\bibitem{BSS84} 
	C. Bardos, R. Santos, and R Sentis, {\em Diffusion approximation and computation of the critical size}, Trans. Amer. Math. Soc. 284, 2, 617-649, 1984
	
	\bibitem{BLP79} 
	A. Bensoussan, J.L. Lions, and G.C. Papanicolaou, {\em Boundary-layers and homogenization of transport processes}, J. Publ. RIMS Kyoto Univ. 15, 53-157, 1979

	\bibitem{Ber}
	F. Berthelin, P. Degond, V. Le Blanc, S. Moutari, J. Royer, M. Rascle, {\em A Traffic-Flow Model with Constraints for the Modeling of Traffic Jams}, Mathematical Models and Methods in Applied Sciences  18, 1269-1298, 2008

	\bibitem{BK18}
	R. Borsche,  A. Klar, {\em A nonlinear discrete velocity relaxation model for traffic flow}, SIAM J. Appl. Math. 78, 5, 2891-2917, 2018
	     
	\bibitem{BK18b}
	R. Borsche,  A. Klar, {\em Kinetic layers and coupling conditions for scalar equations on networks}, Nonlinearity, 31, 7, 3512-3541, 2018	
	
	\bibitem{BK18c}	
	R. Borsche,  A. Klar, {\em Kinetic layers and coupling conditions for macroscopic equations on networks {I}: {T}he wave equation}, SIAM J. Sci. Comput. 40, 3, A1784-A1808, 2018.
		
	\bibitem{BK20}	
	R. Borsche, A. Klar, {\em  A kinetic traffic network  model and its macroscopic limit: diverging  lanes},
	in preparation

	\bibitem{BKKP16} 
	R. Borsche, A. Klar, J.Kall,  T.N.H. Pham, {\em  Kinetic and related macroscopic models for chemotaxis on networks}, M3AS, 26, No. 6, 1219-1242, 2016
	
	\bibitem{BNR14}
	G. Bretti, R. Natalini, M. Ribot, {\em A hyperbolic model of chemotaxis on a network: a numerical study}, ESAIM Math. Model. Numer. Anal., 48(1) ,231--258, 2014.

	\bibitem{CGP05}		
	G.M. Coclite, M.  Garavello,  B. Piccoli, {\em Traffic flow on a road network}, SIAM J. Math. Anal. 36, 6, 1862-1886, 2005.
	
	\bibitem{CG08}
	R.M. Colombo, M. Garavello, {\em On the Cauchy problem for the $p$-system at a junction}, SIAM J. Math. Anal., 39, 1456--1471 2008.

	\bibitem{CM08}
	R.M. Colombo, R.,C. Mauri, {\em Euler system for compressible fluids at a junction}, J. Hyperbolic Differ. Equ., 5(3), 547--568, 2008.

	\bibitem{Cor17}
	A. Corli, L. di Ruvo, L. Malaguti, M. D. Rosini, {\em Traveling waves for degenerate diffusive equations on networks}, NHM 12,3, 339 - 370, 2017.
	
	\bibitem{CGS} 
	F. Coron, F. Golse, C. Sulem, {\em A Classification of Well-posed Kinetic Layer Problems}, CPAM, Vol. 41,  409, 1988.
	
	\bibitem{Dag1} 
	C. F. Daganzo, {\em A behavioral theory of multi-lane traffic flow part I: Long homogeneous freeway
	sections}, Trans. Res. B, 36 (2002), pp. 131–158

	\bibitem{Dag2}
	 C. F. Daganzo, {\em A behavioral theory of multi-lane traffic flow part II: Merges and the onset of congestion}, Trans. Res. B, 36 (2002), pp. 159–169.
	
	\bibitem{G10}
	M. Garavello,
	{\em A review of conservation laws on networks} NHM 5, 3, 565 - 581, 2010  
	
	 \bibitem{GPBook}
	 M. Garavello,  B. Piccoli, {\em Traffic flow on networks}, AIMS Series on Applied Mathematics, 1, American Institute of Mathematical Sciences (AIMS), Springfield, 2006
			
	\bibitem{G08} 
	F. Golse, {\em Analysis of the boundary layer equation in the kinetic theory of gases}, Bull. Inst. Math. Acad. Sin. 3, 1, 211-242, 2008

	 \bibitem{Hel}
	 D.~Helbing,  {\em Gas-kinetic derivation of {N}avier-{S}tokes-like traffic equation}, Physical Review E, 53 (1996), pp.~2366--2381.

	\bibitem{HM09}
	M.~Herty and S.~Moutari, {\em A macro-kinetic hybrid model for traffic flow on road networks}, Comput. Methods Appl. Math., 9, 3,238--252, 2009.
	
	\bibitem{HR}
	M. Herty, M. Rascle, {\em Coupling conditions for a class of second order models for traffic flow}, SIAM Math. Anal. 38,  2, 595–616, 2006.

	\bibitem{HK}
	M. Herty and A. Klar, {\em Modeling, simulation, and optimization of traffic flow networks}, SIAM J. Sci. Comput., 25 (2003), pp. 1066–1087.

	\bibitem{HPRV20}
	M. Herty, G. Puppo, S. Roncoroni, G. Visconti, {\em The BGK approximation of kinetic models for traffic}, Kinetic \& Related Models
 13, 279, 2020.

	\bibitem{HR95}
	H. Holden, N.H. Risebro, {\em A mathematical model of traffic flow on a network of unidirectional roads}, SIAM J. Math. Anal. 26, 4, 999--1017, 1995.
		
	\bibitem{KW97}
	 A.~Klar and R.~Wegener,  {\em Enskog-like kinetic models for vehicular traffic}, J. Stat. Phys., 87 , 91-114, 1997.
  
	\bibitem{L}
	J. P. Lebacque, {\em Les modeles macroscopiques du trafic}, Annales des Ponts, 67 (1993), pp. 28–45.
	
	\bibitem{LS02}
	G. Leugering, Guenter, E.J.P.G. Schmidt, {\em On the modelling and stabilization of flows in networks of open canals}, SIAM J. Control Optim., 41(1), 164--180 2002.

	\bibitem{LX96}
	J.-G. Liu, Z. Xin, {\em Boundary-layer behavior in the fluid-dynamic limit for a nonlinear model Boltzmann Equation}, Arch. Rational Mech. Anal. 135, 61-105, 1996.
	
	\bibitem{NT01}
	R. Natalini, A. Terracina, {\em Convergence of a relaxation approximation to a boundary value problem for conservation laws}, Comm. Partial Differential Equations, 26(7-8), 1235–1252, 2001.

	\bibitem{N99}
	A. Nouri, A. Omrane,J.P. Vila, Boundary conditions for scalar conservation laws from a kinetic point of view. J. Statist. Phys., 94(5-6), 779–804, 1999.

	\bibitem{PSTV17} 
	G. Puppo, M. Semplice, A. Tosin, G. Visconti,{\em Kinetic models for traffic flow resulting in a reduced space of microscopic velocities},Kinet. Relat. Models 10(3), 823--854, 2017.

	\bibitem{PSTV171} 
	G. Puppo, M.Semplice, A. Tosin, G. Visconti, {\em Analysis of multi-population kinetic model for traffic flow}, Commun. Math. Sci. 15(2), 379--412, 2017.

	\bibitem{UTY03} 
	S. Ukai, T. Yang, and S.-H. Yu, {\em Nonlinear boundary layers of the Boltzmann equation. I. Existence}, Comm. Math. Phys. 236,  3, 373-393, 2003.
	
	\bibitem{WX99}
	W.-C. Wang, Z. Xin, {\em Asymptotic limit of initial boundary value problems for conservation laws with relaxational extensions}, Communications on Pure and Applied Mathematics, 51,5  505–535, 1998.
	
	\bibitem{X04}
	W.-Q. Xu, {\em Boundary conditions and boundary layers for a multi-dimensional relaxation model}, Journal of Differential Equations 197,  1, 10, 85-117, 2004.
	
	\bibitem{WY99}
	W.-A. Yong, {\em Boundary conditions for hyperbolic systems with stiff relaxation}, Indiana University Mathematics Journal 48, 1, 115-137,
	1999
	
	
\end{thebibliography}
\end{document}